\documentclass[a4paper,11pt,abstracton]{scrartcl}

\usepackage[T1]{fontenc}             % für "modernere" Schriftkodierung
\usepackage{amsfonts}                % für erweiterte mathematische Schriften
\usepackage{amssymb}                 % für erweiterten Symbolvorrat
\usepackage{textcomp}                % für erweiterten Text-Symbolvorrat
\usepackage{typearea}                % für die Berechnung der Seitenränder
\usepackage{enumerate}               % für komfortablere Aufzählungen
\usepackage{amsmath}
\usepackage{amsthm}
\usepackage{a4}
\usepackage{ifsym}
\usepackage{vmargin}
\usepackage{latexsym}
\usepackage{graphicx}
\usepackage[center]{subfigure}
\usepackage{array}
\usepackage{todonotes}
\usepackage{wrapfig}
\usepackage{url}
\usepackage{scrpage2}
\usepackage{verbatim}
\usepackage{mathrsfs} 							% Calligraphy font
\usepackage[font=small,format=plain,labelfont=bf,up,textfont=it,up]{caption}
\usepackage[numbers]{natbib}
\usepackage{bm} %fette griechische buchstaben mittels \boldsymbol
\usepackage{overpic}
\usepackage{multirow}
\usepackage{fixmath}
\usepackage{epstopdf}
\usepackage{algorithm}
\usepackage[noend]{algpseudocode}
\usepackage{algorithmicx}
\usepackage{xspace}
\usepackage{placeins}
\usepackage{xcolor}
\usepackage{booktabs}
\usepackage{placeins}

\algnewcommand\And{\xspace\textbf{and}\xspace}
\algnewcommand\Or{\xspace\textbf{or}\xspace}
\algnewcommand{\LineComment}[1]{\State \(\triangleright\) #1}

\graphicspath{{fig/}}

\newcommand{\N}{\mathbb{N}}

\newcommand{\Sp}{\mathbb{S}}

\newcommand{\R}{\mathbb{R}}

\newcommand{\curlG}{\overrightarrow{\textrm{curl}}_{\Gamma}}

\newcommand{\e}{\operatorname{e}}
\newcommand{\dif}{{\ensuremath{\,\mathrm{d}}}}

\numberwithin{equation}{section}
\numberwithin{theorem}{section}

\numberwithin{figure}{section}

\begin{document}

\title{Efficient Solution of Time-Domain Boundary Integral Equations Arising in Sound-Hard Scattering}
\author{
A. Veit\thanks{aveit@uchicago.edu, Department of Computer Science, The University of Chicago, 1100 E. 58th Street, Chicago, IL 60637 US} 
\and 
M. Merta\thanks{michal.merta@vsb.cz, V\v{S}B-Technical University of Ostrava, 
17. listopadu 2172/15, 708 33 Ostrava, Czech Republic}
\and 
J. Zapletal\thanks{(jan.zapletal@vsb.cz), V\v{S}B-Technical University of Ostrava, 
17. listopadu 2172/15, 708 33 Ostrava, Czech Republic}
\and 
D. Luk\'{a}\v{s}\thanks{(dalibor.lukas@vsb.cz), V\v{S}B-Technical University of Ostrava, 
17. listopadu 2172/15, 708 33 Ostrava, Czech Republic}
}
\date{}
\maketitle

\begin{abstract}
We consider the efficient numerical solution of the three-dimensional wave equation with Neumann boundary conditions via time-domain boundary integral equations. A space-time Galerkin method with $C^\infty$-smooth, compactly supported basis functions in time and piecewise polynomial basis functions in space is employed. We discuss the structure of the system matrix and its efficient parallel assembly. Different preconditioning strategies for the solution of the arising systems with block Hessenberg matrices are proposed and investigated numerically. Furthermore, a C++ implementation parallelized by OpenMP and MPI in shared and distributed memory, respectively, is presented. The code is part of the boundary element library BEM4I. Results of numerical experiments including convergence and scalability tests up to a thousand cores on a cluster are provided. The presented implementation shows good parallel scalability of the system matrix assembly. Moreover, the proposed algebraic preconditioner in combination with the FGMRES solver leads to a significant reduction of the computational time.\vspace{\baselineskip}


\end{abstract}

\section{Introduction}
We are concerned with the efficient numerical solution of time-dependent scattering phenomena in unbounded domains. Specifically, we consider the time-dependent, three-dimensional wave equation in the case of a sound-hard scatterer modelled by Neumann boundary conditions. We formulate and solve the arising problem in terms of time-domain boundary integral equations (TDBIEs). The main advantage of this approach is that the problem (originally posed in a three-dimensional unbounded domain) is reduced to the two-dimensional (bounded) surface of the scatterer. Different discretization techniques have been introduced and intensively studied to efficiently solve TDBIEs. These methods are typically based on the Galerkin discretization in space and utilize either collocation, convolution quadrature or a Galerkin scheme for the time discretization. We refer to \cite{costabel2004time,BanjaiSchanz,Geranmayeh,HaDuong,DuongLudwigTerrasse} and the references therein for an overview of existing methods and to \cite{Sayas} for a detailed introduction to the theory of time-domain boundary integral equations.

Here, we will employ the space-time Galerkin method to solve the arising equations numerically. This approach originated from the groundbreaking work of Bamberger and Ha Duong \cite{BamDuongSoft,BamDuong}. They introduced coercive space-time variational formulations for acoustically soft and hard scatterers and showed stability and convergence of the resulting Galerkin schemes with piecewise polynomial ansatz spaces. The main difficulty in this approach is the accurate computation of the matrix entries. If standard piecewise polynomial basis functions are used, special quadrature techniques taking into account the complicated shape of the arising 4-dimensional integration domains are necessary. To circumvent this problem, we use $C^\infty$-smooth and compactly supported basis functions for the time-discretization. These were introduced in \cite{Sauter5, SauterVeitIMA} for the Dirichlet problem and it was shown that this choice allows an accurate approximation of the matrix entries using standard quadrature schemes. As the consequence of the simplified computation of the system matrix the use of nonequidistant timesteps and higher order approximation spaces in time became feasible \cite{SauterVeit2013}.

In this paper we employ the same type of ansatz functions for the Neumann problem using the variational formulation derived in \cite{BamDuong}. We focus on equidistant timesteps and investigate the implications on the structure of the system matrix. Due to the overlap of the temporal basis functions the resulting linear system that needs to be solved admits a block Hessenberg structure. We compare a conventional GMRES solver with a GMRES solver preconditioned by deflations (see \cite{ErhBurPoh96}) and show numerically that the necessary number of iterations is significantly reduced. In Section \ref{Sec:AlgPrecond} we furthermore introduce an experimental algebraic preconditioner that exploits the block Hessenberg structure of the system matrix. We perform various numerical experiments that show the performance of this preconditioner in combination with the FGMRES method developed in \cite{saad1993}.

We present an efficient parallel implementation of the above-mentioned discretization and solution strategies in Section \ref{Sec:Implementation}. The code is a part of the boundary element library BEM4I \cite{BEM4I}. It is based on C++ and uses hybrid parallelization by OpenMP and MPI to accelerate the evaluation of the discretized space-time boundary integral operators. The implementation leverages the repeating pattern of the system matrices to minimize computational cost as well as memory requirements. We briefly describe the structure of the BEM4I library and  provide details about the parallel system matrix assembly in Section~\ref{Sec:parallel_assembly}. Results of numerical experiments demonstrating scalability of the computation in a distributed memory architecture are provided in Section \ref{Sec:NumExp}. 

\section{Integral Formulation of the Wave Equation}
Let $\Omega\subset\mathbb{R}^{3}$ be a Lipschitz domain with the boundary denoted by $\Gamma
$. We consider the homogeneous wave equation
\begin{subequations}
\label{FullProblem}
\end{subequations}
\begin{equation}
\partial_{t}^{2}u-\Delta u=0\hspace{3mm}\text{in }\Omega\times\left[
0,T\right]  \tag{\ref{FullProblem}a}\label{WaveEquation}
\end{equation}
with homogeneous initial conditions
\begin{equation}
u(\cdot,0)=\partial_{t}u(\cdot,0)=0\hspace{3mm}\text{in }\Omega\tag{\ref{FullProblem}b}\label{InitialConditions}
\end{equation}
and Neumann boundary conditions
\begin{equation}
\partial_nu:=\frac{\partial u}{\partial n}=g\hspace{3mm}\text{on }\Gamma\times\left[  0,T\right]  \tag{\ref{FullProblem}c}\label{BoundaryConditions}
\end{equation}
on a time interval $\left[  0,T\right]  $ for $T>0$, where $n$ denotes the unit outward normal vector. In applications, $\Omega$
is often the unbounded exterior of a bounded domain. For such problems, the
method of boundary integral equations is an elegant tool where the partial
differential equation is transformed to an equation on the bounded surface
$\Gamma$. 

We employ an ansatz as a \textit{double layer potential} for the
solution $u$,
\begin{align}
u(x,t):=&D\phi(x,t):=-\frac{1}{4\pi}\int_{\Gamma}\frac{n_y\cdot (x-y)}{\| x-y\|}\left(\frac{\phi(y,t-\|x-y\|)}{\|x-y\|^2}+\right.\notag\\
&\left.\frac{\partial_{t}\phi(y,t-\|x-y\|)}{\|x-y\|} \right)\dif\Gamma_y,\hspace{3mm}\ (x,t)\in\Omega\backslash\Gamma\times\left[  0,T\right]
\label{RetardedPotential}%
\end{align}
with an unknown density function $\phi$. $D$ is also referred to as
\textit{retarded double layer potential} due to the retarded time argument
$t-\Vert x-y\Vert$ connecting the time and space variables.

The ansatz (\ref{RetardedPotential}) satisfies the wave equation
(\ref{WaveEquation}) and the initial conditions \eqref{InitialConditions}. Therefore, the unknown density function $\phi$ has to be determined such that the Neumann boundary conditions \eqref{BoundaryConditions} are satisfied.
For this the normal derivative of the double layer potential has to be extended to the boundary $\Gamma$ which can be done continuously for sufficiently smooth functions across smooth points of $\Gamma$. We therefore define the hypersingular operator
\begin{equation}
\label{Hypersingular}
W v(x,t) := \lim_{x^{+}\in\Omega\rightarrow x} n_x \cdot \nabla_{x^{+}}D v(x^{+},t)
\end{equation}
for $(x,t)\in\Gamma\times[0,T]$, where the limit is taken in the sense of distributions. \begin{comment} We remark that, for sufficiently smooth functions $v$, the normal derivative of the double layer potential is continuous when $x\in\Omega\backslash\Gamma$ approaches the boundary $\Gamma$, i.e.,\todo{Shouldn't there be something like "lim from interior equals lim form exterior"?}
\[
 \lim_{x^{+}\in\Omega\rightarrow x}\partial_n D v (x^{+},t) = W v (x,t).
\]
\end{comment}
In order to find the unknown density function $\phi$ in \eqref{RetardedPotential} such that \eqref{BoundaryConditions} is satisfied we thus consider the boundary integral equation
\begin{equation}
\label{BIE}
W\phi = g \quad \text{on }\Gamma\times [0,T]. 
\end{equation}
To solve this equation numerically we introduce a weak formulation of \eqref{BIE} following \cite{BamDuong}. A suitable space-time variational formulation is given by: Find $\phi$ in a Sobolev space $V$ such that
\begin{align}
a(\phi,\zeta)&:=\int_{0}^{T}\int_{\Gamma}\int_{\Gamma}\left\lbrace \frac{n_x\cdot n_y}{4\pi\Vert x-y\Vert}\ddot{\phi}(y,t-\Vert x-y\Vert)\dot{\zeta}(x,t) \right.\notag\\
&\qquad\left.+\frac{\curlG\phi(y,t-\Vert x-y\Vert)\cdot \curlG\dot{\zeta}(x,t) }{4\pi\Vert x-y\Vert} \right\rbrace \dif\Gamma_{y}\dif\Gamma_{x}\dif t\nonumber\\
&=\int_{0}^{T}\int_{\Gamma}g(x,t)\dot{\zeta}(x,t)\dif\Gamma_{x}\dif t=:b(\zeta) \label{VarForm}%
\end{align}
for all $\zeta\in V$, where we denote by $\dot{\phi}$ and $\ddot{\phi}$  the first and second derivatives with respect to time. Here, $\curlG\phi$ is the tangential rotation of the function $\phi$ defined as
\[
\curlG\phi(x) := n_x \times \nabla \tilde{\phi}(x),
\]
where $\tilde{\phi}$ is defined in a tubular neighbourhood of $\Gamma$ by
\[
\tilde{\phi}(x+\varepsilon n_x):=\phi(x)
\]
for $x\in\Gamma$ (see also \cite{BamDuong, nedelec2001acoustic}).

\section{Numerical Discretization}
We discretize the variational problem \eqref{VarForm} using a Galerkin method in space and time.  Therefore, we replace the infinite dimensional space $V$ by a finite dimensional subspace $V_{\operatorname{Galerkin}}$ spanned by $L$ basis functions $\{b_{i}\}_{i=1}^{L}$ in time and $M$
basis functions $\{\varphi_{j}\}_{j=1}^{M}$ in space. This leads to the
discrete ansatz
\begin{equation}
\phi_{\operatorname{Galerkin}}(x,t)=\sum_{i=1}^{L}\sum_{j=1}^{M}\alpha_{i}%
^{j}\varphi_{j}(x)b_{i}(t),\hspace{3mm}(x,t)\in\Gamma\times\left[  0,T\right]
,\label{DiscreteAnsatz}%
\end{equation}
with the unknown coefficients $\alpha_{i}^{j}$. Plugging \eqref{DiscreteAnsatz} into the variational formulation \eqref{VarForm} and using the basis functions $b_k$ and $\varphi_l$ as test functions leads to the linear system
\begin{equation}
\underline{\underline{\mathbf{A}}}\cdot\underline{\boldsymbol{\alpha}}= \underline{\mathbf{g}},
\label{LinearSystemOriginal}
\end{equation}
where the block matrix $\underline{\underline{\mathbf{A}}}\in\R^{LM\times LM}$, the unknown coefficient vector $\underline{\boldsymbol{\alpha}}\in\R^{LM}$ and the right-hand side vector $\underline{\mathbf{g}}\in\R^{LM}$ can be partitioned according to
\begin{equation}
\underline{\underline{\mathbf{A}}}:=%
\begin{bmatrix}
\mathbf{A}_{1,1} & \mathbf{A}_{1,2} & \cdots & \mathbf{A}_{1,L}\\
\mathbf{A}_{2,1} & \mathbf{A}_{2,2} & \cdots & \mathbf{A}_{2,L}\\
\vdots & \vdots & \ddots & \vdots\\
\mathbf{A}_{L,1} & \mathbf{A}_{L,2} & \cdots & \mathbf{A}_{L,L}%
\end{bmatrix}
,\qquad\underline{\boldsymbol{\alpha}}:=%
\begin{bmatrix}
\boldsymbol{\alpha}_{1}\\
\boldsymbol{\alpha}_{2}\\
\vdots\\
\boldsymbol{\alpha}_{L}%
\end{bmatrix}
,\qquad\underline{\mathbf{g}}:=%
\begin{bmatrix}
\mathbf{g}_{1}\\
\mathbf{g}_{2}\\
\vdots\\
\mathbf{g}_{L}%
\end{bmatrix}
,\label{A}%
\end{equation}
with
\[
\mathbf{A}_{k,i}\in\mathbb{R}^{M\times M},\quad \boldsymbol{\alpha}_i\in\R^{M} ,\quad\mathbf{g}_{k}\in\mathbb{R}%
^{M}\quad\text{ for }i,k\in\{1,\cdots,L\}.
\]
Individual entries are given by
\begin{align}
\label{integrals}
\mathbf{A}_{k,i}(j,l)= &\int_{0}^{T}\int_{\Gamma}\int_{\Gamma}\left\lbrace \frac{n_x\cdot n_y}{4\pi\Vert x-y\Vert}\,\varphi _{j}(y)\ddot{b}_i(t-\Vert x-y\Vert)\,\varphi_{l}(x)\dot{b}_k(t) \right.\notag\\
&\qquad\left.+\frac{\curlG\varphi _{j}(y) \cdot \curlG\varphi_{l}(x)  }{4\pi\Vert x-y\Vert}b_i(t-\Vert x-y\Vert)\dot{b}_k(t) \right\rbrace \dif\Gamma_{y}\dif\Gamma_{x}\dif t
\end{align}
and
\begin{equation*}
\boldsymbol{\alpha}_i(j)=\left(\alpha_i^j\right)_{j=1}^M,\quad \mathbf{g}_{k}(l)=\int_{0}^{T}\int_{\Gamma }g(x,t)\,\varphi_{l}(x)\,\dot{b}_{k}(t)\dif\Gamma _{x}\dif t,
\end{equation*}
respectively. We rewrite \eqref{integrals} by introducing univariate functions  
\begin{equation}
\psi_{k,i}(r):=\int_{0}^{T}\ddot{b}_{i}(t-r)\dot{b}_{k}(t)\dif t,\qquad \tilde{\psi}_{k,i}(r):=\int_{0}^{T}{b}_{i}(t-r)\dot{b}_{k}(t)\dif t
\label{psis}
\end{equation}
(see \cite{SauterVeit2013}) and obtain
\begin{align}
\mathbf{A}_{k,i}(j,l)= &\int_{\Gamma}\int_{\Gamma}\left\lbrace \frac{n_x\cdot n_y}{4\pi\Vert x-y\Vert}\,\varphi _{j}(y)\,\varphi_{l}(x) \psi_{k,i}(\|x-y\|) \right.\notag\\
&\qquad\left.+\frac{\curlG\varphi _{j}(y) \cdot \curlG\varphi_{l}(x)  }{4\pi\Vert x-y\Vert} \tilde{\psi}_{k,i}(\|x-y\|) \right\rbrace \dif\Gamma_{y}\dif\Gamma_{x}.
\label{MatrixEntries}
\end{align}

The accurate computation of the matrix entries \eqref{MatrixEntries} is problematic in the space-time Galerkin approach. If piecewise polynomial basis functions in time are used as proposed in \cite{BamDuong}, the integrand in \eqref{MatrixEntries} is only piecewise smooth which makes standard quadrature techniques prohibitively expensive. In \cite{Sauter5}, $C^\infty$-smooth and compactly supported temporal shape functions $b_{i}$ were proposed. It could be shown (see \cite{SauterVeitIMA,SauterVeit2013,KhSaVe}) that this choice significantly simplifies the accurate approximation of integrals as in \eqref{MatrixEntries} and as a consequence allows the use of nonequidistant stepsizes as well as higher-order approximation spaces in time. Since our goal in this paper is a fast solver, we restrict ourselves to equidistant timesteps as the computational complexity and the memory requirements are significantly lower in this case. We denote the timesteps by 
\[
t_i := i\cdot\Delta t,\quad \text{with }\Delta t:=\frac{T}{N-1},\, i=0,\ldots N-1,
\]
where $N$ is the number of timesteps. In the following we briefly recall the definition of the temporal basis functions for the special case of equidistant timesteps. We define
\[
f\left(  t\right)  :=\left\{
\begin{array}
[c]{ll}%
\frac{1}{2}\operatorname{erf}\left(  2\operatorname{artanh}t\right)  +\frac
{1}{2} & \left\vert t\right\vert <1,\\
0 & t\leq-1,\\
1 & t\geq1
\end{array}
\right.
\]
and note that $f\in C^{\infty}\left(  \mathbb{R}\right)  $. We scale and shift $f$ in order to obtain a (left) cutoff function
\[
f_{i}\left(  t\right)  :=f\left(  2\frac{t-t_i}{\Delta t}-1\right),\quad\text{where }f_{i}\left(  t\right)=\left\{
\begin{array}
[c]{ll}%
0 & t\leq t_i,\\
1 & t\geq t_{i+1}.
\end{array}
\right.
\]
We obtain a bump function on the interval $\left[ t_{i-1},t_{i+1}\right]  $ with midpoint
$t_i$ by
\[
\rho_{i}\left(  t\right)  :=\left\{
\begin{array}
[c]{ll}%
f_{i-1}\left(  t\right)   & t_{i-1}\leq t\leq t_{i},\\
1-f_{i}\left(  t\right)   & t_{i}\leq t\leq t_{i+1},\\
0 & \text{otherwise.}%
\end{array}
\right.
\]
A smooth partition of unity of the interval $[0,T]$ is then defined by
\[
\mu_{1}:=1-f_{0},\quad\mu_{N}:=f_{N-2},\quad\forall i\in \{2, \ldots, N-1\}:\mu_{i}:=\rho_{i-1}.
\]
Smooth and compactly supported basis functions in time can now be
obtained by multiplying these partition of unity functions with suitably
scaled Legendre polynomials $P_m$ of degree $m$ (see \cite{Sauter5} for details):
\begin{alignat}{2}
&b_{1,m}:=8\mu_{1}(t)\,\left(\frac{t}{\Delta t}\right)^{2}P_{m}\left(  \frac{2}{\Delta t
}t-1\right)\quad  && m=0,\ldots,p,\nonumber\\
&b_{i,m}:=\mu_{i}(t)P_{m}\left(  \frac{t-t_{i-2}}%
{\Delta t}-1\right)  && m=0,\ldots,p,\hspace{1mm}i=2,\ldots,N-1,\label{BasisFunctions}\\
&b_{N,m}:=\mu_{N}(t)P_{m}\left(  2\frac{t-t_{N-2}%
}{\Delta t}-1\right) \quad && m=0,\ldots,p, \nonumber
\end{alignat}
where $p$ controls the order of the method in time. We will use the above basis functions in time for the Galerkin approximation \eqref{DiscreteAnsatz}. Note that the definition is slightly different to the one in \cite{Sauter5}. Here we simply use $p+1$ basis functions in the first interval. This choice leads to a better asymptotic convergence rate of the method in the interval $[t_0,t_1]$  since we have the same number of basis functions as for the other time intervals but additionally use the a priori knowledge about the solution $u(\cdot,0)=\partial_{t}u(\cdot,0)=0$, which is reflected in the factor $t^2$ in the basis functions (see \cite{Babu}). This choice will simplify the implementation. In order to use the functions \eqref{BasisFunctions} as basis functions in the discrete ansatz \eqref{DiscreteAnsatz} a suitable numeration has to be introduced. Here we use 
\[
b_i := b_{\left\lceil \frac{i}{p+1} \right\rceil, \operatorname{mod}(i-1,p+1)}\quad \text{for} \quad i=1,\ldots ,L=N\cdot p.
\]
For the discretization in space we use standard piecewise polynomial (typically linear) basis functions $\varphi_{j}$.

\subsection{Efficient representation and evaluation of $\psi_{k,i}$ and $\tilde{\psi}_{k,i}$}
\label{Sec:PsiEvaluation}

An efficient handling of the functions $\psi_{k,i}$ and $\tilde{\psi}_{k,i}$ in \eqref{psis} is crucial for a successful implementation of the algorithm since they have to be evaluated numerous times during the approximation of the matrix entries \eqref{MatrixEntries}. In \cite{SauterVeit2013} we propose to approximate this type of functions for each $k,i\in\{1,\ldots,L\}$ with piecewise Chebyshev polynomials. This approximation is efficient due to the smoothness of $\psi_{k,i}$ and $\tilde{\psi}_{k,i}$ and it can be evaluated efficiently with Clenshaw's recurrence formula~\cite{Clenshaw}. Here we furthermore exploit the fact that we only use constant timesteps. The number of approximations that have to be precomputed in this case is only $\mathcal{O}(p^2)$ compared to $\mathcal{O}(N^2p^2)$ for variable stepsizes in time. 

Let the indices $i,k\in\{1,\ldots,L\}$ be arbitrary but fixed and let $\tilde{i},\tilde{k}\in\{1,\ldots,N\}$,  $m_1,m_2\in\{0,\ldots,p\}$ be such that
\begin{equation}
b_i\equiv b_{\tilde{i},m_1},\quad b_k\equiv b_{\tilde{k},m_2}.
\label{basisEquiv}
\end{equation}
We first consider the case where $2\leq \tilde{i},\tilde{k} \leq N-1$, i.e., basis functions that are associated with ``inner'' timesteps. Then simple calculus shows that
\begin{align}
\tilde{\psi}_{k,i}(r) &= \int_{t_{\tilde{k}-2}}^{t_{\tilde{k}}}{b}_{\tilde{i},m_1}(t-r)\dot{b}_{\tilde{k},m_2}(t)\dif t\nonumber\\
&=\int_{t_{\tilde{k}-2}}^{t_{\tilde{k}}}{b}_{2,m_1}(t-t_{\tilde{i}-2}-r)\dot{b}_{2,m_2}(t-t_{\tilde{k}-2})\dif t\nonumber\\
&=:\tilde{\xi}_{m_1,m_2}(r+t_{\tilde{i}-2}-t_{\tilde{k}-2}),\label{psiTildeRed}
\end{align}
where
\[
\tilde{\xi}_{m_1,m_2}(\alpha) :=\int_{0}^{2\Delta t}{b}_{2,m_1}(t-\alpha)\dot{b}_{2,m_2}(t)\dif t.
\]
Thus the task of approximating $\tilde{\psi}_{k,i}$ for $p+2\leq k,i \leq L-p-2$ is reduced to approximating $\tilde{\xi}_{m_1,m_2}$ for $0\leq m_1,m_2 \leq p$ on its support $[-2\Delta t,2\Delta t]$ and evaluating these functions according to \eqref{psiTildeRed}. Completely similarly we obtain that
\[
\psi_{k,i}(r) = \xi_{m_1,m_2}(r+t_{\tilde{i}-2}-t_{\tilde{k}-2}),\quad\text{with}\quad\xi_{m_1,m_2}(\alpha) :=\int_{0}^{2\Delta t}\ddot{b}_{2,m_1}(t-\alpha)\dot{b}_{2,m_2}(t)\dif t.
\]

\begin{figure}[!ht]
\centering
\subfigure[$m_1=0,m_2=0$]{
\centering
\includegraphics[width=0.47\textwidth]{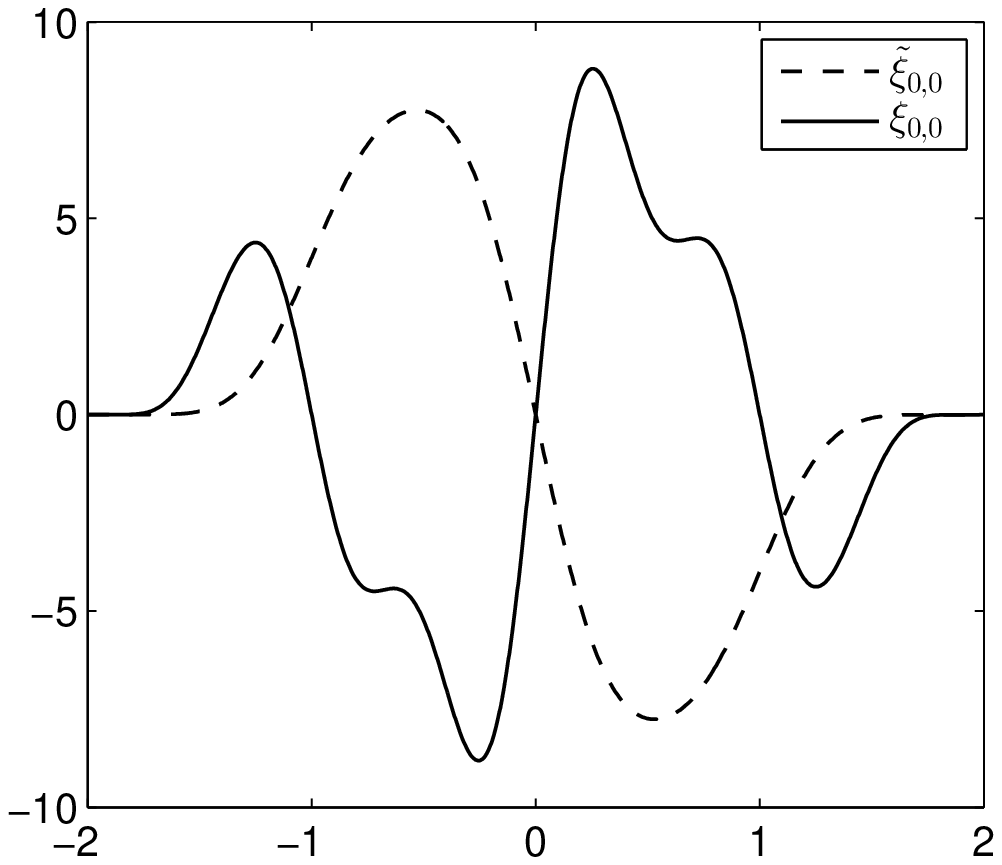}
} \hfil
\subfigure[$m_1=1,m_2=0$]{
\centering
\includegraphics[width=0.46\textwidth]{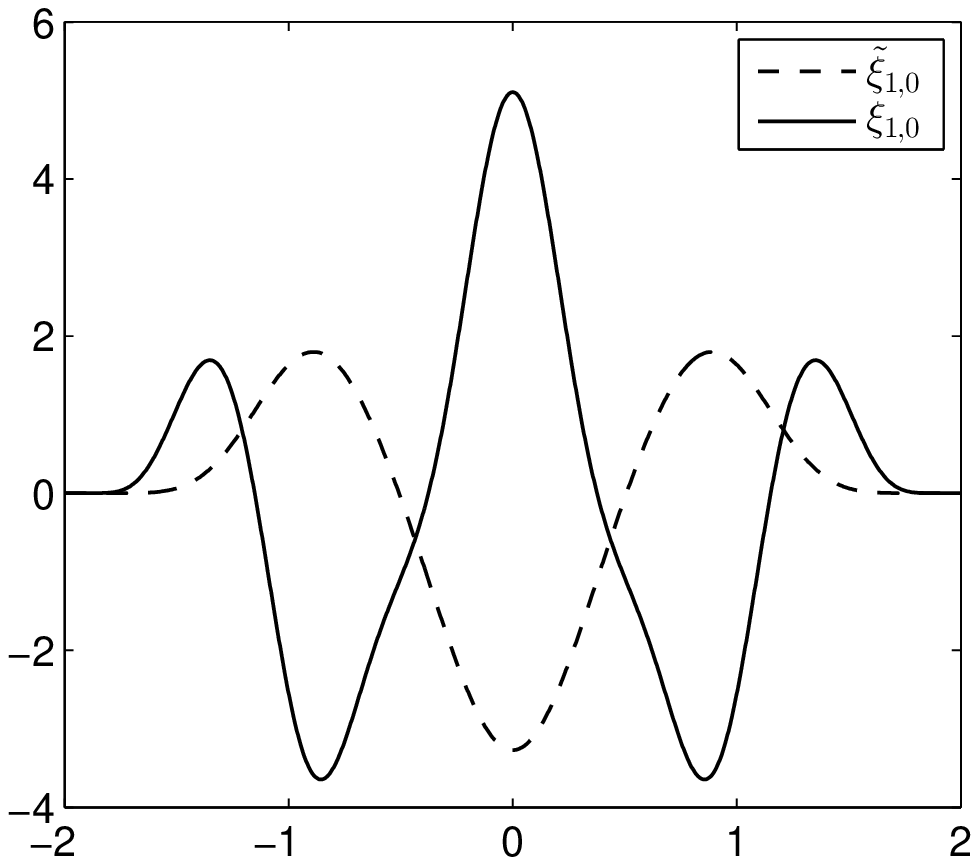}
} \caption{$\xi_{m_1,m_2}$ and $\tilde{\xi}_{m_1,m_2}$ for $\Delta t = 1$. For visualization purposes $\tilde{\xi}_{0,0}$ and $\tilde{\xi}_{1,0}$ were multiplied by a factor 8.}%
\label{Fig:psis}
\end{figure}

Figure \ref{Fig:psis} illustrates the prototype functions $\tilde{\xi}_{m_1,m_2}$ and $\xi_{m_1,m_2}$ for $\Delta t=1$. It suggests that these functions are either even or odd depending on $m_1$ and $m_2$. Indeed it is easy to see that
\begin{equation}
\tilde{\xi}_{m_1,m_2}(\alpha) = (-1)^{m_1+m_2+1}\tilde{\xi}_{m_1,m_2}(-\alpha)
\label{symmetryXi}
\end{equation}
due to the symmetry of the basis functions $b_{2,m}(t)$ with respect to $t=\Delta t$. The same formula holds for $\xi_{m_1,m_2}$. As a consequence, approximations of $\xi_{m_1,m_2}$ and $\tilde{\xi}_{m_1,m_2}$ (e.g., with piecewise Chebyshev polynomials) have only to be computed and stored for $\alpha\in [0,2\Delta t]$ since
\begin{equation}
\tilde{\psi}_{k,i}(r) =  \operatorname{sign}\left(r+t_{\tilde{i}-2}-t_{\tilde{k}-2}\right)^{m_1+m_2+1}\tilde{\xi}_{m_1,m_2}\left(\left|r+t_{\tilde{i}-2}-t_{\tilde{k}-2}\right|\right),
\label{FormulaForPsi}
\end{equation}
where $k,\tilde{k},m_2$ and $i,\tilde{i},m_1$ are related via \eqref{basisEquiv}. An analogous formula holds for $\psi_{k,i}$.

Lastly we want to point out a symmetry with respect to $m_1$ and $m_2$. Partial integration yields
\begin{equation}
\tilde{\xi}_{m_1,m_2}(-\alpha) = -\tilde{\xi}_{m_2,m_1}(\alpha)\quad\text{and}\quad \xi_{m_1,m_2}(-\alpha) = -\xi_{m_2,m_1}(\alpha).
\label{symmetryXip}
\end{equation}
Therefore, only $\frac{1}{2}(p+1)(p+2)$ functions $\tilde{\xi}_{m_1,m_2}$ and $\xi_{m_1,m_2}$ have actually to be approximated. More importantly, this relation has an impact on the structure of the system matrix $\underline{\underline{\mathbf{A}}}$ as we will show in Section \ref{Sec:MatrixStructure}.

So far we only considered the case $2\leq \tilde{i},\tilde{k}\leq N-1$. If the basis functions are associated with the first or the last timestep, similar prototype functions can be defined. Since this is analogous to the procedure described above we omit the details here.

\subsection{Structure of the system matrix $\underline{\underline{\mathbf{A}}}$}
\label{Sec:MatrixStructure}
The solution of \eqref{VarForm} using the discrete ansatz
\eqref{DiscreteAnsatz} leads to a linear system with $L\cdot M$ unknowns. The special choice of basis functions in time as well as the use of equidistant timesteps has several implications on the structure of the system matrix $\underline{\underline{\mathbf{A}}}$. Throughout this section we assume again that $k,\tilde{k},m_2$ and $i,\tilde{i},m_1$ are related via \eqref{basisEquiv}. We denote the matrix block $\mathbf{A}_{k,i}$ from \eqref{integrals} by $\mathbf{A}^{m_2,m_1}_{\tilde{k},\tilde{i}}$ to highlight the dependence on the timestep and the order of the involved basis functions in time. Furthermore, we define the matrix block
\begin{equation}
\mathbf{\tilde{A}}_{\tilde{k},\tilde{i}}:=%
\begin{bmatrix}
\mathbf{A}^{0,0}_{\tilde{k},\tilde{i}}  & \cdots & \mathbf{A}^{0,p}_{\tilde{k},\tilde{i}}\\
\vdots  & \ddots & \vdots\\
 \mathbf{A}^{p,0}_{\tilde{k},\tilde{i}}  & \cdots &  \mathbf{A}^{p,p}_{\tilde{k},\tilde{i}}
\end{bmatrix}\in\R^{(p+1)M\times (p+1)M}.
\end{equation}
We first remark that
\[
\operatorname{supp}\psi_{k,i} = \operatorname{supp}\tilde{\psi}_{k,i} =  [t_{\tilde{k}-2} - t_{\tilde{i}}, t_{\tilde{k}} - t_{\tilde{i}-2}]\cap \R_{\geq 0},
\]
where we formally set $t_{-1}=0$ and $t_{N}=T$. Thus
\[
\psi_{k,i} \equiv \tilde{\psi}_{k,i} \equiv  0 \text{ on }\R\quad \text{ for }t_{\tilde{k}} \leq t_{\tilde{i}-2}.
\]
This shows that the resulting system matrix $\underline{\underline{\mathbf{A}}}\in\R^{LM\times LM}$ admits the block Hessenberg structure
\begin{equation}
\underline{\underline{\mathbf{A}}}=%
\begin{bmatrix}
\mathbf{\tilde{A}}_{1,1} & \mathbf{\tilde{A}}_{1,2} & \mathbf{0} & \mathbf{0} & \cdots & \mathbf{0}\\
\mathbf{\tilde{A}}_{2,1} & \mathbf{\tilde{A}}_{2,2} & \mathbf{\tilde{A}}_{2,3} & \mathbf{0} & \cdots & \mathbf{0}\\
\vdots & \vdots& \vdots& \ddots & \ddots & \vdots\\
\vdots & \vdots& \vdots& \vdots & \ddots & \mathbf{0}\\
\mathbf{\tilde{A}}_{N-1,1} & \mathbf{\tilde{A}}_{N-1,2} & \cdots& \cdots & \cdots & \mathbf{\tilde{A}}_{N-1,N}\\
\mathbf{\tilde{A}}_{N,1} & \mathbf{\tilde{A}}_{N,2} & \cdots& \cdots & \cdots & \mathbf{\tilde{A}}_{N,N}
\end{bmatrix},
\label{Apartitioned}
\end{equation}
where $\mathbf{0}\in \R^{(p+1)M\times (p+1)M}$ denotes the zero matrix. An important consequence of \eqref{FormulaForPsi} is that for $2\leq \tilde{i},\tilde{k} \leq N-1$ the functions $\psi_{k,i}$ and $\tilde{\psi}_{k,i}$ only depend on the difference $t_{\tilde{i}-2}-t_{\tilde{k}-2}$ and therefore
\begin{equation}
\mathbf{\tilde{A}}_{\tilde{k},\tilde{i}} = \mathbf{\tilde{A}}_{\tilde{k}-\tilde{i}+2,2}\text{ for }\tilde{k}-\tilde{i}\geq 0 \quad\text{and}\quad \mathbf{\tilde{A}}_{\tilde{k},\tilde{i}} = \mathbf{\tilde{A}}_{2,3}\text{ for }\tilde{k}-\tilde{i}= -1.
\label{relationAtildeBlock}
\end{equation}
Thus only $\mathcal{O}(N)$ matrix blocks $\mathbf{\tilde{A}}_{\tilde{k},\tilde{i}}$ have actually to be computed and stored. In order to further reduce the complexity we remark that for $2\leq \tilde{i},\tilde{k} \leq N-1$ the formula \eqref{symmetryXip} together with \eqref{symmetryXi} shows that
\begin{equation}
\mathbf{A}^{m_1,m_2}_{\tilde{k},\tilde{i}} = (-1)^{m_1+m_2} \mathbf{A}^{m_2,m_1}_{\tilde{k},\tilde{i}}.
\label{relationAblock}
\end{equation}
In order to represent $\mathbf{\tilde{A}}_{\tilde{k},\tilde{i}}$ we thus only have to compute and store  $\frac{1}{2}(p+1)(p+2)$ matrix blocks $\mathbf{A}^{m_2,m_1}_{\tilde{k},\tilde{i}}$. Finally we remark that the blocks $\mathbf{A}^{m_2,m_1}_{\tilde{k},\tilde{i}}$ themselves are sparse (see also \cite{Sauter5,SauterVeitIMA}) since
\begin{equation}
\mathbf{A}^{m_2,m_1}_{\tilde{k},\tilde{i}}(j,l) = 0\quad\text{if}\quad \operatorname{supp}\psi_{k,i}\cap [\operatorname{mindist}_{j,l},\operatorname{maxdist}_{j,l}]=\emptyset,
\label{SparsitySmallBlock}
\end{equation}
where
\begin{align*}
\operatorname{mindist}_{j,l}&:=\min\left\lbrace \|x-y\|,\,x\in\operatorname{supp}\varphi_l,\,y\in\operatorname{supp}\varphi_j\right\rbrace\\
\operatorname{maxdist}_{j,l}&:=\max\left\lbrace \|x-y\|,\,x\in\operatorname{supp}\varphi_l,\,y\in\operatorname{supp}\varphi_j\right\rbrace.
\end{align*}
Note that in the case $T>\operatorname{diam}\left(\Gamma\right)$, \eqref{SparsitySmallBlock} implies that $\mathbf{A}^{m_2,m_1}_{\tilde{k},\tilde{i}}=\mathbf{0}$ and therefore $\mathbf{\tilde{A}}_{\tilde{k},\tilde{i}}=\mathbf{0}$ if $t_{\tilde{k}-2}-t_{\tilde{i}}>\operatorname{diam}\left(\Gamma\right)$.

It is clear that a computational and memory efficient implementation should avoid the explicit construction of $\underline{\underline{\mathbf{A}}}$. Instead, only those matrix blocks $\mathbf{A}^{m_2,m_1}_{\tilde{k},\tilde{i}}$ should be computed that are necessary to recover arbitrary entries of $\underline{\underline{\mathbf{A}}}$ via \eqref{relationAtildeBlock} and \eqref{relationAblock}. Furthermore, these blocks have to be stored in a suitable sparse matrix storage format. Since the solution of the arising linear system \eqref{LinearSystemOriginal} is typically computed using iterative solvers it is also necessary to develop a routine for the efficient multiplication of  $\underline{\underline{\mathbf{A}}}$ with a vector that incorporates the special structure of $\underline{\underline{\mathbf{A}}}$.

%\begin{comment}
%In order to evaluate the matrix entries \eqref{MatrixEntries} efficiently by standard quadrature rules the functions $\psi_{k,i}$ and $\tilde{\psi}_{k,i}$ have to be evaluated several times. This is expensive since these functions are defined themselves by an integral which has to be approximated. However due to the smoothness of the basis functions in time, $\psi_{k,i}$ and $\tilde{\psi}_{k,i}$ are also smooth functions on the whole real line. We therefore approximate these functions for every $k,i\in\lbrace 1,\ldots ,L \rbrace$ on their support accurately with a polynomial, which can then be evaluated efficiently. \todo{We should discuss this point in more detail - cite Sauter, Veit: Adaptive Time Discretization for Retarded Potentials (Michal)?}
%\end{comment}

\section{Solution of the linear system}
\label{Sec:SolLinSys}
In contrast to classical schemes using piecewise polynomial temporal basis functions the resulting system matrix is not lower triangular with Toeplitz structure and FFT-type methods cannot be used for the solution of the linear system. Therefore an efficient iterative solver has to be utilized and an appropriate preconditioner of the system has to be developed. The GMRES algorithm is one of the possible choices \cite{Saad1986gmres}. However, since the system matrix $\underline{\underline{\mathbf{A}}}$ is global in time and space and therefore of large dimension, the restarted version of the algorithm GMRES($m$) is usually necessary to keep the computational and memory requirements reasonable. To speed up the convergence we present two possible preconditioning techniques based on deflation and a recursive algebraic approach.

\subsection{Restarted GMRES preconditioned by deflation}

It is known that restarts of the algorithm lead to a slower convergence than in the case of the full-GMRES due to the loss of information on the smallest Ritz values \cite{huang1989some}. The restarted GMRES preconditioned by deflation (DGMRES($m,l$)) aims to keep this information by approximating the invariant subspace corresponding to the smallest eigenvalues of the system matrix. After each restart the invariant subspace is increased by $l$ new vectors. It has been observed that DGMRES($m,l$) can restore the convergence even in cases where the original restarted algorithm stalls \cite{ErhBurPoh96}. 

Let $|\lambda_1| \leq |\lambda_2|\leq\ldots\leq|\lambda_{LM}|$ be the eigenvalues of $\underline{\underline{\mathbf{A}}}$. The desired preconditioner has the form 
\[\widehat{\mathbf{M}}:=\mathbf{I}_{LM}+\mathbf{U}(1/|\lambda_{LM}|\mathbf{T}-\mathbf{I}_r)\mathbf{U}^T,\]
where $\mathbf{T}:=\mathbf{U}^T\underline{\underline{\mathbf{A}}}\mathbf{U}$, $\mathbf{I}_{LM},\,\mathbf{I}_r$ are identity matrices of appropriate dimensions, and $\mathbf{U}$ is an orthonormal basis of the invariant subspace corresponding to the smallest $r$ eigenvalues of $\underline{\underline{\mathbf{A}}}$. It can be proven that 
\[\widehat{\mathbf{M}}^{-1}=\mathbf{I}_{LM}+\mathbf{U}(|\lambda_{LM}|\mathbf{T}^{-1}-\mathbf{I}_r)\mathbf{U}^T\]
and the eigenvalues of $\underline{\underline{\mathbf{A}}}\widehat{\mathbf{M}}^{-1}$ are $\lambda_{r+1}, \lambda_{r+2}, \ldots, \lambda_{LM}, |\lambda_{LM}|$. Thus, the smallest $r$ eigenvalues of $\underline{\underline{\mathbf{A}}}$ are removed and replaced by $|\lambda_{LM}|$ with multiplicity $r$.

Since it is inefficient to assemble the preconditioner $\widehat{\mathbf{M}}$ accurately, we aim to set up its suitable approximation $\mathbf{M}$. The GMRES algorithm produces the basis $\mathbf{V}_k$ of the current Krylov subspace and the Hessenberg matrix $\mathbf{H}_k=\mathbf{V}_k^T\underline{\underline{\mathbf{A}}}\mathbf{V}_k$ representing the restriction of $\underline{\underline{\mathbf{A}}}$ onto this subspace. In order to construct $\mathbf{M}$, the Schur decomposition of $\mathbf{H}_k$ is performed. Using the Schur vectors $\mathbf{S}$ corresponding to the smallest eigenvalues of $\mathbf{H}_k$ we approximate the Schur vectors of $\underline{\underline{\mathbf{A}}}$ as $\mathbf{U}\approx\mathbf{V}_k\mathbf{S}$. Moreover, the largest eigenvalue of $\mathbf{H}_k$ is used to approximate the largest eigenvalue of~$\underline{\underline{\mathbf{A}}}$.

Using this procedure the algorithm increases $\mathbf{U}$ by $l$ vectors after each restart until the maximal dimension $r$ of the invariant subspace is achieved. Since the matrices $\mathbf{H}_k$ are dense and of relatively small dimensions, the Schur decomposition is done efficiently using the BLAS routines. Moreover, it is not necessary to explicitly assemble the matrix $\mathbf{M}^{-1}$ and its actions are carried out as a sequence of dense matrix-vector multiplications. For more detailed description we refer the reader to \cite{ErhBurPoh96}.

\subsection{An algebraic preconditioner for block Hessenberg systems}
\label{Sec:AlgPrecond}
In this section we propose a preconditioner for the linear system \eqref{LinearSystemOriginal} that makes use of the block Hessenberg structure of $\underline{\underline{\mathbf{A}}}$. Let
\begin{equation}
\underline{\underline{\mathbf{A}}}\underline{\underline{\mathbf{M}}}^{-1}(\underline{\underline{\mathbf{M}}}\cdot\underline{\boldsymbol\alpha}) =\underline{\mathbf{g}} 
\label{precondSystem}
\end{equation}
be the preconditioned system where $\underline{\underline{\mathbf{M}}}$ denotes the preconditioner. We assume that $\underline{\underline{\mathbf{A}}}$ is partitioned according to \eqref{Apartitioned}. We choose $\underline{\underline{\mathbf{M}}}$ to be the matrix that coincides with $\underline{\underline{\mathbf{A}}}$ except that the matrix block $\tilde{\mathbf{A}}_{\lceil N/2\rceil,\lceil N/2\rceil+1}$ is replaced by a zero matrix of the corresponding size, i.e.,
\[
\mathbf{\underline{\underline{M}}}:=\begin{bmatrix}
 \mathbf{\underline{\underline{M}}}_{1,1} & \mathbf{\underline{\underline{0}}} \\
 \mathbf{\underline{\underline{M}}}_{2,1}  &  \mathbf{\underline{\underline{M}}}_{2,2} 
\end{bmatrix},
\]
where
\begin{align*}
\mathbf{\underline{\underline{M}}}_{1,1}&:=\left(\tilde{\mathbf{A}}_{\tilde{k},\tilde{i}}\right)_{\tilde{k},\tilde{i}=1\ldots\lceil N/2\rceil},\\
\mathbf{\underline{\underline{M}}}_{2,1}&:=\left(\tilde{\mathbf{A}}_{\tilde{k},\tilde{i}}\right)_{\tilde{k}=\lceil N/2+1\rceil\ldots N ,\, \tilde{i}:=1\ldots \lceil N/2\rceil},\\ \mathbf{\underline{\underline{M}}}_{2,2}&=\left(\tilde{\mathbf{A}}_{\tilde{k},\tilde{i}}\right)_{\tilde{k},\tilde{i}=\lceil N/2+1\rceil\ldots N}.
\end{align*}
This choice of $\underline{\underline{\mathbf{M}}}$ as a preconditioner in \eqref{precondSystem} is motivated by the Woodbury matrix identity which states that the inverse of a rank-$k$ perturbed matrix can be computed by doing a rank-$k$ correction to the inverse of the original matrix. Since $\underline{\underline{\mathbf{M}}}$ is a rank-$(p+1)M$ perturbation of $\underline{\underline{\mathbf{A}}}$ it follows that $\underline{\underline{\mathbf{A}}}\underline{\underline{\mathbf{M}}}^{-1}$ is a rank-$(p+1)M$ perturbation of the identity matrix. This suggests that the application of an iterative solver like GMRES is more efficient for the preconditioned system \eqref{precondSystem} than for the original one.

In order to apply the preconditioner within an iterative solver, $\underline{\underline{\mathbf{M}}}^{-1}$ is not needed explicitly. It is however crucial that the actions of $\underline{\underline{\mathbf{M}}}^{-1}$ on a vector $\underline{\mathbf{r}}=\left(\underline{\mathbf{r}}_1^T,\underline{\mathbf{r}}_2^T\right)^T$ can be computed efficiently. Due to the block triangular structure of $\underline{\underline{\mathbf{M}}}$ it is easy to see that
\[
\mathbf{\underline{\underline{M}}}^{-1}\mathbf{\underline{r}} = \begin{bmatrix}
 \mathbf{\underline{\underline{M}}}_{1,1}^{-1}\mathbf{\underline{r}}_1  \\
 \mathbf{\underline{\underline{M}}}_{2,2}^{-1}\left(\mathbf{\underline{r}}_2-\mathbf{\underline{\underline{M}}}_{2,1} \mathbf{\underline{\underline{M}}}_{1,1}^{-1}\mathbf{\underline{r}}_1\right)  
\end{bmatrix}
\]
Therefore, the problem of computing $\mathbf{\underline{\underline{M}}}^{-1}\mathbf{\underline{r}} $ boils down to the problem of evaluating  $ \mathbf{\underline{\underline{M}}}_{1,1}^{-1}\mathbf{\underline{v}}$ and $ \mathbf{\underline{\underline{M}}}_{2,2}^{-1}\mathbf{\underline{w}}$ for some vectors $\mathbf{\underline{v}}$ and $\mathbf{\underline{w}}$. This is equivalent to the solution of  linear systems of the form
\begin{equation}
\mathbf{\underline{\underline{M}}}_{1,1}\mathbf{\underline{x}}=\mathbf{\underline{v}}\quad \text{ and}\quad  \mathbf{\underline{\underline{M}}}_{2,2}\mathbf{\underline{x}}=\mathbf{\underline{w}}. 
\label{innerSystems}
\end{equation}
Note that $\mathbf{\underline{\underline{M}}}_{1,1}$ and $\mathbf{\underline{\underline{M}}}_{2,2}$ are again block Hessenberg matrices with basically half the size of the original matrix $\mathbf{\underline{\underline{A}}}$. Thus, these smaller systems can again be solved iteratively using a preconditioner of the same form. This strategy can be applied recursively until the matrices that have to be inverted are just the diagonal blocks $\mathbf{\tilde{A}}_{k,k}$. On the lowest level the resulting linear systems then should be solved either exactly or with a standard iterative solver.

\begin{algorithm}[ht]
\footnotesize
\begin{algorithmic}
\Require System matrix $\underline{\underline{\mathbf{A}}}$, RHS vector $\underline{\mathbf{g}}$, relative precision $\varepsilon$, max. number of iterations $m_{\mathrm{it}}$
\State Set max. recursion depth of the preconditioner $m_{\mathrm{r}}$
\State Set max. number of inner iterations $m_{\mathrm{it}}^{\mathrm{in}}$
\State Initialize the current recursion level $r:=1$
\State\Call{AssemblePreconditioner}{$\underline{\underline{\mathbf{A}}},\underline{\underline{\mathbf{M}}}$}
\State \Call{FGMRES}{$\underline{\underline{\mathbf{A}}},\underline{\boldsymbol\alpha}, \underline{\mathbf{g}}, \varepsilon, m_{\mathrm{it}},\underline{\underline{\mathbf{M}}}$}
\Ensure Solution vector $\underline{\boldsymbol\alpha}$
\\\hrulefill
\Function{AssemblePreconditioner}{$\mathbf{A}$, $\mathbf{M}$}
\LineComment{assembles the preconditioner $\mathbf{M}$ for the block Hessenberg matrix $\mathbf{A}$}
\State $\mathbf{M}_{1,1}:=(\mathbf{A}_{k,i})_{k,i=1,\ldots,\lceil N/2\rceil}$
\State $\mathbf{M}_{2,1}:=(\mathbf{A}_{k,i})_{k=\lceil N/2+1\rceil,\ldots,N,i=1,\ldots\lceil N/2\rceil}$
\State $\mathbf{M}_{2,2}:=(\mathbf{A}_{k,i})_{k,i=\lceil N/2+1\rceil,\ldots,N}$
\State $\mathbf{M}:=\begin{pmatrix}
 \mathbf{M}_{1,1} & \mathbf{0} \\
 \mathbf{M}_{2,1}  &  \mathbf{M}_{2,2} 
\end{pmatrix}$
\EndFunction
\\\hrulefill
\Function{ApplyPreconditioner}{$\mathbf{M}$, $\mathbf{x}$, $\mathbf{y}$}
\LineComment{applies the preconditioner $\mathbf{M}$ on a vector $\mathbf{x}$ and stores the result in $\mathbf{y}$}
\State $\mathbf{M}_{1,1}:=(\mathbf{M}_{k,i})_{k,i=1,\ldots,\lceil N/2\rceil}$
\State $\mathbf{M}_{2,1}:=(\mathbf{M}_{k,i})_{k=\lceil N/2+1\rceil,\ldots,N,i=1,\ldots\lceil N/2\rceil}$
\State $\mathbf{M}_{2,2}:=(\mathbf{M}_{k,i})_{k,i=\lceil N/2+1\rceil,\ldots,N}$
\State $\mathbf{x}=[\mathbf{x}_1, \mathbf{x}_2]^T, \mathbf{y}=[\mathbf{y}_1, \mathbf{y}_2]^T$
\If{$r\leq m_{\mathrm{r}}$}
\State\Call{AssemblePreconditioner}{$\mathbf{M}_{1,1}$, $\mathbf{M}_{1}$}
\State\Call{AssemblePreconditioner}{$\mathbf{M}_{2,2}$, $\mathbf{M}_{2}$}
\Else
\State $\mathbf{M}_{1}=\mathbf{M}_{2}:=\mathbf{I}$
\EndIf
\State $r:=r+1$
\State \Call{FGMRES}{$\mathbf{M}_{1,1}$, $\mathbf{y}_1$, $\mathbf{x}_1$, $\varepsilon$, $m_{\mathrm{it}}^{\mathrm{in}}$, $\mathbf{M}_{1}$} 
\State \Call{FGMRES}{$\mathbf{M}_{2,2}$, $\mathbf{y}_2$, $\mathbf{x}_2-\mathbf{M}_{2,1}\mathbf{y}_1$, $\varepsilon$, $m_{\mathrm{it}}^{\mathrm{in}}$, $\mathbf{M}_{2}$} 
\State return $\mathbf{y}=[\mathbf{y}_1,\mathbf{y}_2]^T$
\EndFunction
\\\hrulefill
\Function{FGMRES}{$\mathbf{A}$, $\mathbf{x}$, $\mathbf{y}$, $\varepsilon$, $m_{\mathrm{it}}$, $\mathbf{M}$}
\LineComment{\parbox[t]{0.8\textwidth}{solves the system $\mathbf{A}\mathbf{x}=\mathbf{y}$ using FGMRES algorithm with relative precision $\varepsilon$, maximum number of iterations $m_{\mathrm{it}}$, and preconditioner $\mathbf{M}$}}
\State $\vdots$
\EndFunction
\end{algorithmic}
\caption{Solution of the system using the recursive preconditioner}
\label{alg:preconditioner}
\end{algorithm}

If the inner systems \eqref{innerSystems} and subsequent smaller systems in the recursive process are solved very accurately the application of this preconditioner is typically too expensive and the resulting solver is too time consuming. Thus we suggest to only employ a limited number of iterations of an inner solver to approximate the solution of \eqref{innerSystems}. Since the classical GMRES algorithm does not allow the preconditioner to change at every iteration (and thus forbids the usage of an inexact solver inside the preconditioner) we use a flexible inner-outer preconditioned variant of GMRES (FGMRES) presented in \cite{saad1993} to solve both \eqref{LinearSystemOriginal} and the systems \eqref{innerSystems}. FGMRES allows a variable preconditioner at the cost of double memory requirements but with the same arithmetic complexity. Numerical experiments in Section~\ref{Sec:NumExp} indicate that only a few iterations of the inner FGMRES solver during the application of $\mathbf{\underline{\underline{M}}}^{-1}$ are necessary to significantly reduce the number of outer iterations and the solution time.

A simplified description of the solver using the recursive preconditioner is listed in Algorithm~\ref{alg:preconditioner}. First, the maximum recursion level of the preconditioner and the number of iterations of the inner FGMRES solver are set. Next, the preconditioner is assembled and the FGMRES is called to solve the outer system~\eqref{LinearSystemOriginal}. During its application the preconditioner itself employs $m_{\mathrm{it}}^{\mathrm{in}}$ iterations of the preconditioned FGMRES to approximate the solution of the systems \eqref{innerSystems}. This repeats recursively until the maximum level of recursion is reached. On the lowest level the inner systems are solved with FGMRES without a preconditioner.

The application of this recursive preconditioner is experimental. Although the numerical benchmarks in Section \ref{Sec:NumExp} show promising results a rigorous convergence analysis is necessary. Moreover, only a sequential version has been tested and its proper parallelization is yet to be performed. 

\section{Implementation}
\label{Sec:Implementation}
A parallel solver for sound scattering problems based on the approach described in the previous sections has been implemented in the BEM4I library \cite{BEM4I}. The code is written in C++ in an object oriented way and utilizes OpenMP and MPI for parallelization in shared and distributed memory, respectively. All classes are templated to support various indexing and scalar types.

\begin{figure}[htb]
\begin{center}
\includegraphics[width=0.9\textwidth]{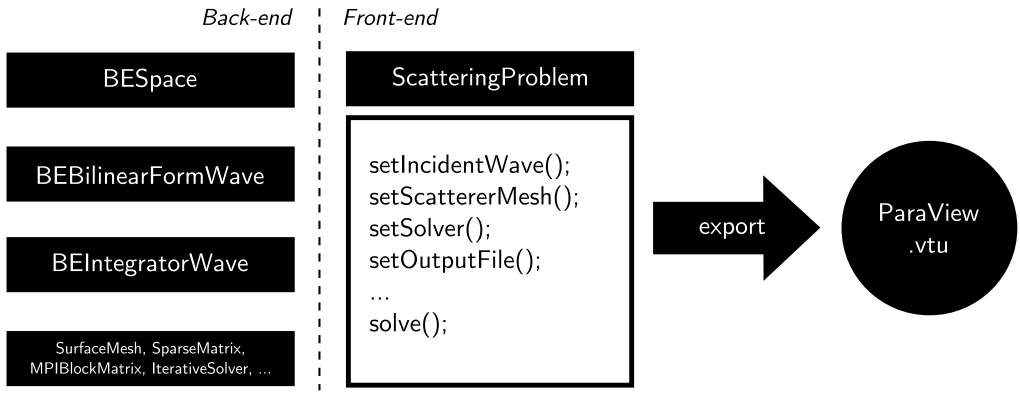}
\end{center}
\caption{Structure of the solver for wave scattering problems in the BEM4I library.}
\label{fig:LibraryStructure}
\end{figure}

The structure of the solver is depicted in Figure \ref{fig:LibraryStructure}. Its core consists of a set of classes responsible for the assembly of system matrices. 

\begin{enumerate}
\item The \texttt{BESpace} class holds the information about the order of spatial and temporal basis and test functions. It also stores the object of the underlying surface mesh.
\item The purpose of the \texttt{BEBilinearFormWave} class is to assemble appropriate system matrices. Several matrix types are supported, including non-distributed sparse matrices (suitable for small problems solved by a direct solver), non-distributed block sparse matrices, and block sparse matrices distributed among computational nodes using MPI. The usage of block system matrices reduces memory and computational requirements since it does not duplicate assembly and storage of identical blocks. The assembly of individual blocks is parallelized by OpenMP.
\item The \texttt{BEIntegratorWave} class is responsible for the assembly of local system matrices, i.e., the quadrature over pairs of elements. It takes care of evaluation of smooth temporal basis and testing functions and their Chebyshev approximations.
\end{enumerate}
Besides these main classes the solver utilizes classes representing surface meshes, full, sparse, and block matrices, direct and iterative solvers (including GMRES, DGMRES, and FGMRES), classes for the evaluation of potential operators, etc. A user solves a problem either by direct manipulation of the above-mentioned classes or using the interface provided in the \texttt{ScatteringProblem} class. The solution is exported to the ParaView \texttt{.vtu} file format.

\subsection{Assembly of distributed system matrix}
\label{Sec:parallel_assembly}
Special care has to be taken during the parallel assembly of the system matrix to exploit its structure and minimize memory and computational requirements. One has to take into account its properties described in Section~\ref{Sec:MatrixStructure}, i.e., the block Hessenberg format~\eqref{Apartitioned}, duplication of blocks following from \eqref{relationAtildeBlock}, sparsity and the special structure of individual blocks \eqref{relationAblock}, as well as the fact that $\mathbf{\tilde{A}}_{\tilde{k}, \tilde{i}}=\mathbf{0}$ for $t_{\tilde{k}-2,\tilde{i}}-t_{\tilde{i}} > \operatorname{diam}\left(\Gamma\right)$. In the following we describe the assembly of a distributed block system matrix. Similar principles can be applied to assemble a non-distributed block matrix suitable for problems of smaller dimension.

\begin{figure}[ht]
\begin{center}
\includegraphics[width=0.95\textwidth]{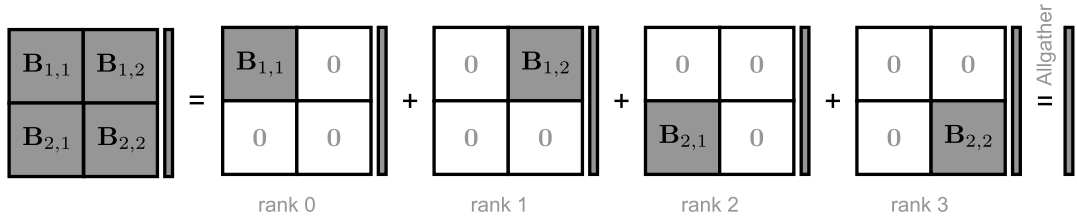}
\end{center}
\caption{Multiplication by a distributed matrix \texttt{MPIBlockMatrix}.}
\label{fig:distributedApply}
\end{figure}

To represent a system matrix we use the \texttt{MPIBlockMatrix} class of the BEM4I library. The class serves as a simple distributed container for non-distributed matrices, such as \texttt{FullMatrix}, \texttt{SparseMatrix}, or any linear operator implementing the \texttt{Apply} method on a vector. When an instance of the class is created in the MPI environment, each MPI process is only assigned a subset of matrix blocks. The matrix-vector multiplication is performed in parallel and the local results are gathered, summed, and distributed among all processes. The multiplied vector is replicated on every process (see Figure~\ref{fig:distributedApply}).

\begin{figure}[t]
\begin{center}
\includegraphics[width=0.85\textwidth]{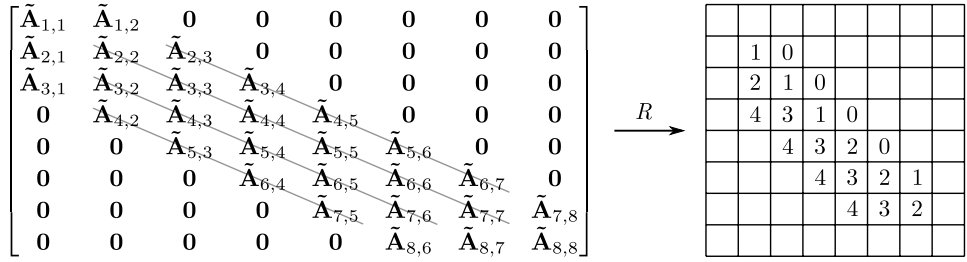}
\end{center}
\caption{Example of a system matrix and distribution of its inner blocks $\mathbf{\tilde{A}}_{2,2}, \ldots, \mathbf{\tilde{A}}_{N-1,N-1}$ among five processes with ranks $\left\{0, 1, \ldots, 4\right\}$. Identical blocks are connected by a line (here $N=8, \tilde{n}_z=3$).}
\label{fig:distribution}
\end{figure}

In the first step of the matrix assembly, the distribution of non-zero blocks $\mathbf{\tilde{A}}_{\tilde{k},\tilde{i}}$ among $P$ available MPI processes is determined. Assuming as above that the time interval $[0, T]$ is discretized into equidistant timesteps of the length $\Delta t:=T/(N-1)$, the number of inner non-zero blocks $\mathbf{\tilde{A}}_{\tilde{k},\tilde{i}}, \tilde{k},\tilde{i}\in\{2, \ldots, N-1\}$ can be obtained as
\begin{equation}
%n_z:=-\frac{\tilde{n}_z^2}{2}+\tilde{n}_z N+\frac{\tilde{n}_z}{2}+N-1,
%n_z:=-\frac{\tilde{n}_z^2}{2}+\tilde{n}_z N-\frac{3}{2}\tilde{n}_z+N-3,
%n_z:=\frac{1}{2}N(N-1)-\frac{N-\tilde{n}_z-2}{2}\max\{N-\tilde{n}_z-1,0\},
n_z:=\frac{1}{2}(N^2-N-4)-\frac{N-\tilde{n}_z-2}{2}\max\{N-\tilde{n}_z-1,0\},
\end{equation}
where 
\begin{equation}
\tilde{n}_z:=\min\left\{N,\left\lceil\frac{2\cdot\text{diam}(\Gamma)}{\Delta t}\right\rceil+2\right\}
\end{equation}
is the number of non-zero blocks in the first matrix column. Each MPI process is assigned approximately $n_z/P$ inner non-zero blocks. The exact distribution of blocks among processes is defined by a mapping $R: \N\times\N\rightarrow\N_0$ assigning the MPI rank of the owner to every block. Two factors have to be considered when defining $R$ -- proper load balance during the matrix-vector multiplication and duplication of blocks $\mathbf{\tilde{A}}_{\tilde{k},2}, \tilde{k}\in\{2, \ldots, N-2\}$ and $\mathbf{\tilde{A}}_{2,3}$ due to the relation \eqref{relationAtildeBlock}. Therefore, the distribution is performed diagonal-wise in order to store identical blocks on the same process in one memory location. However, some blocks are duplicated among several processes to ensure balanced load during the matrix-vector multiplication. An example of the distribution of inner blocks among processes is depicted in Figure~\ref{fig:distribution}. After the assignment of inner blocks, the blocks in the first and last rows and columns are mapped to processes with the smallest number of blocks assigned. Note that this distribution does not take into account differences in fill-in of various blocks.

After the mapping $R$ is defined, the blocks $\mathbf{\tilde{A}}_{\tilde{k}, \tilde{i}}$ are assembled. Due to the relation \eqref{relationAtildeBlock} the total number of blocks to be assembled is at most $3N$ (the first two columns, the last row, and the blocks $\mathbf{\tilde{A}}_{2,3}$ and $\mathbf{\tilde{A}}_{N-1,N}$), therefore the complexity is reduced from $\mathcal{O}(N^2)$ to $\mathcal{O}(N)$. Each unique non-zero block is assembled and distributed in the following way:

\begin{enumerate}
\item An instance of the \texttt{SparseMatrix} class representing the current block $\mathbf{\tilde{A}}_{\tilde{k}, \tilde{l}}$ is created. 
\item Approximations of temporal basis functions associated with the current block are computed using the relations from \cite{SauterVeit2013} and Section~\ref{Sec:PsiEvaluation}.
\item Using \eqref{SparsitySmallBlock} the non-zero pattern of the current block $\mathbf{\tilde{A}}_{k,l}$ is precomputed together with pairs of mesh elements contributing to the non-zero entries. This computation is parallelized by MPI and OpenMP.
\item The pairs of elements are evenly split among available MPI processes.
\item Every MPI process iterates through its element pairs and through the order of the temporal basis functions and computes its local contribution to the block $\mathbf{\tilde{A}}_{\tilde{k}, \tilde{i}}$ using \eqref{MatrixEntries}. In order to further save computational time and memory property \eqref{relationAblock} can be employed. The iteration through element pairs is parallelized in shared memory by OpenMP.
\item Each process sends its partial results to the owner(s) of the block determined by the mapping $R$. The owner combines the results to form the current block. If any identical block is owned by the process, only a reference to the assembled block is copied to the appropriate memory location and the data is not duplicated.
\item All other processes delete their partial results.
\end{enumerate}

\section{Numerical experiments}
\label{Sec:NumExp}
In the following section we present results of numerical experiments performed on the Anselm cluster located at the IT4Innovations National Supercomputing Centre, Czech Republic. The cluster consists of 188 compute nodes, each of them equipped with two 8-core Intel Xeon E5-2665 processors and 64 GB of RAM. The~theoretical peak performance of the cluster is 82 Tflop/s.

\subsection{Convergence tests}
In this section we test the convergence of the space-time Galerkin method introduced above. In \cite{Veit_diss}, analytic solutions of direct and indirect formulations for Dirichlet and Neumann problems on the unit sphere were derived. Such solutions turned out to be very useful for numerical experiments and the validation of the algorithm and also serve here as reference solutions. We recall these solutions for the problem \eqref{BIE} for the sake of completeness.

Let us assume that the scatterer is the unit sphere, i.e., $\Gamma = \mathbb{S}^2$, and that the right-hand side in \eqref{BIE} is of the form
\begin{equation}
g(x,t)=g(t)Y_n^m(x) \quad \text{for }(x,t)\in\mathbb{S}^2\times [0,T],
\label{eq:rhs}
\end{equation}
where $Y_n^m(x)$ are spherical harmonics of degree $n$ and order $m$. We assume $g(t)$ to be causal, i.e., $g(t)=0$ for $t\leq 0$ and that $\dot{g} (0) = 0$. It can be shown (see \cite{nedelec2001acoustic}) that $Y_n^m(x)$ are eigenfunctions of the hypersingular operator in the frequency domain, i.e.,
\[
\mathcal{L}(W)(s)Y_n^m(x) = \lambda_n(s)Y_n^m(x),
\]
where $\mathcal{L}$ denotes the Laplace transform. The function $\lambda_n(s)$ can be expressed explicitely in terms of spherical Bessel functions and spherical Hankel functions of the first kind. Assuming that the solution of \eqref{BIE} for the special right-hand side \eqref{eq:rhs} also admits the form $\phi(t)Y_n^m(x)$ and exploiting the fact that $W$ is a convolution with respect to the time variable we obtain
\begin{equation}
\phi(t) = \int_0^t g(\tau)\mathcal{L}^{-1}\left\lbrace\frac{1}{\lambda_n}\right\rbrace(t-\tau) \dif\tau.
\label{formula_phi}
\end{equation}
Evaluating the inverse Laplace transform in the formula above leads to explicit formulae for $\phi$. In the case $n=0$, i.e., the right-hand side $g(x,t)=g(t)$ is purely time-dependent, the solution of \eqref{BIE} is purely time-dependent as well and is given by
\begin{align*}
\phi(x,t)&=-2\int_{0}^{t}g(t-\tau) \cosh (\tau)\dif\tau \\
&\quad+2\sum_{k=1}^{\lfloor t/2\rfloor} \sum_{l=1}^{k} (-1)^{k+1}\int_{2k}^{t}  c_{k,l}\,(\tau-2k)^{k-l+1} \e^{\tau-2k}   \dot{g}(t-\tau) \dif\tau,
\end{align*}
where
\[
c_{k,l}:=\binom{k-1}{l-1}\frac{2^{k-l}}{(k-l+1)!}.
\]
In the case $n=1$, i.e., the right-hand side is of the form $g(x,t)=g(t)Y_1^m(x)$  the solution of \eqref{BIE} is given by
\[
\phi(x,t)=\left(-2\int_{0}^{t}g(t-\tau)\cosh(\tau)\cos(\tau)\dif\tau\right) Y_1^m(x)
\]
for $t\in [0,2[$. Formulae for large $n$ and $t$ are typically complicated and furthermore difficult to evaluate due to numerical cancellation. In this case the one-dimensional problem
\[
 \int_0^t \mathcal{L}^{-1}\left\lbrace\lambda_n\right\rbrace(t-\tau)\phi(\tau) \dif\tau = g(t)
 \]
can be accurately and efficiently solved using the convolution quadrature in order to obtain approximations of \eqref{formula_phi} (see, e.g., \cite{BanjaiSchanz}).\vspace{\baselineskip}

\begin{figure}[htb]
\begin{center}
\includegraphics[width=0.5\textwidth]{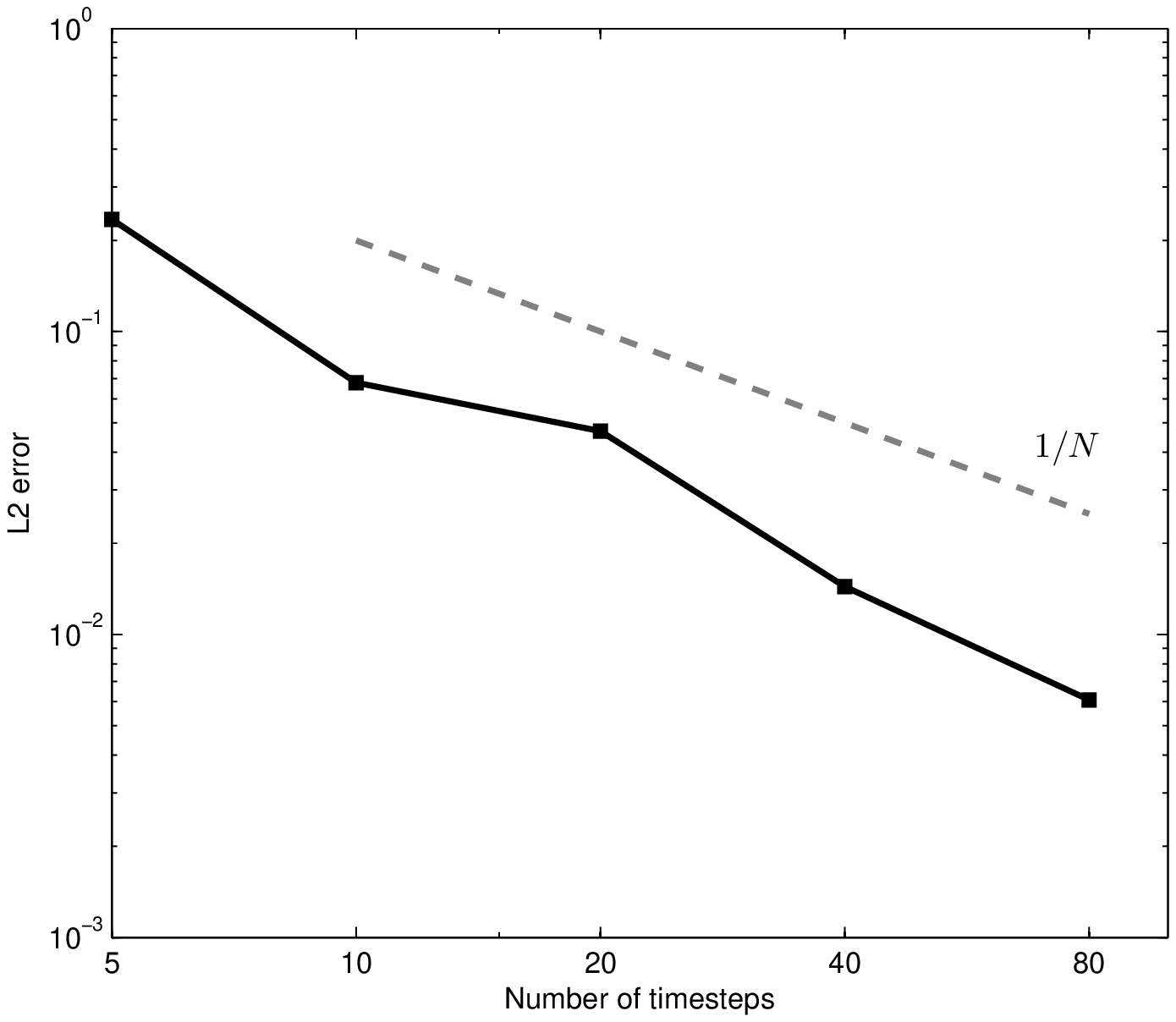}\includegraphics[width=0.5\textwidth]{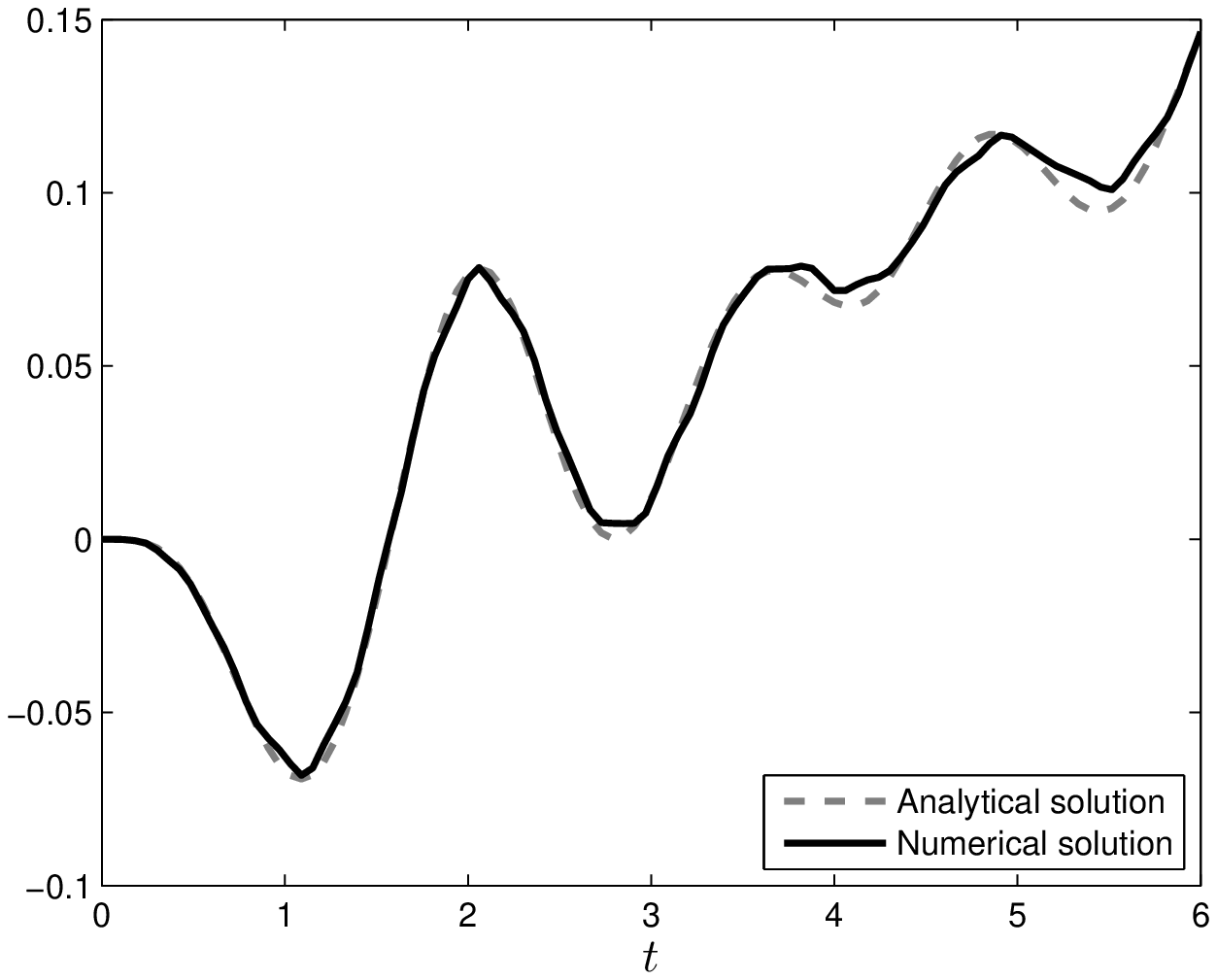}
\end{center}
\caption{Left: log-log scale plot of $\|\Phi-\Phi_{\mathrm{Galerkin}}\|_{L^2(\Gamma;\,[0,T])}$ for $g(t)=\sin(3t)t^2\e^{-t}, n=0, T=6$, and local polynomial approximation space of degree $p=1$. Right: numerical solution in $x\in\Gamma$ obtained using 20 timesteps and $p=1$.}
\label{fig:error1}
\end{figure}

\begin{figure}[htb]
\begin{center}
\includegraphics[width=0.5\textwidth]{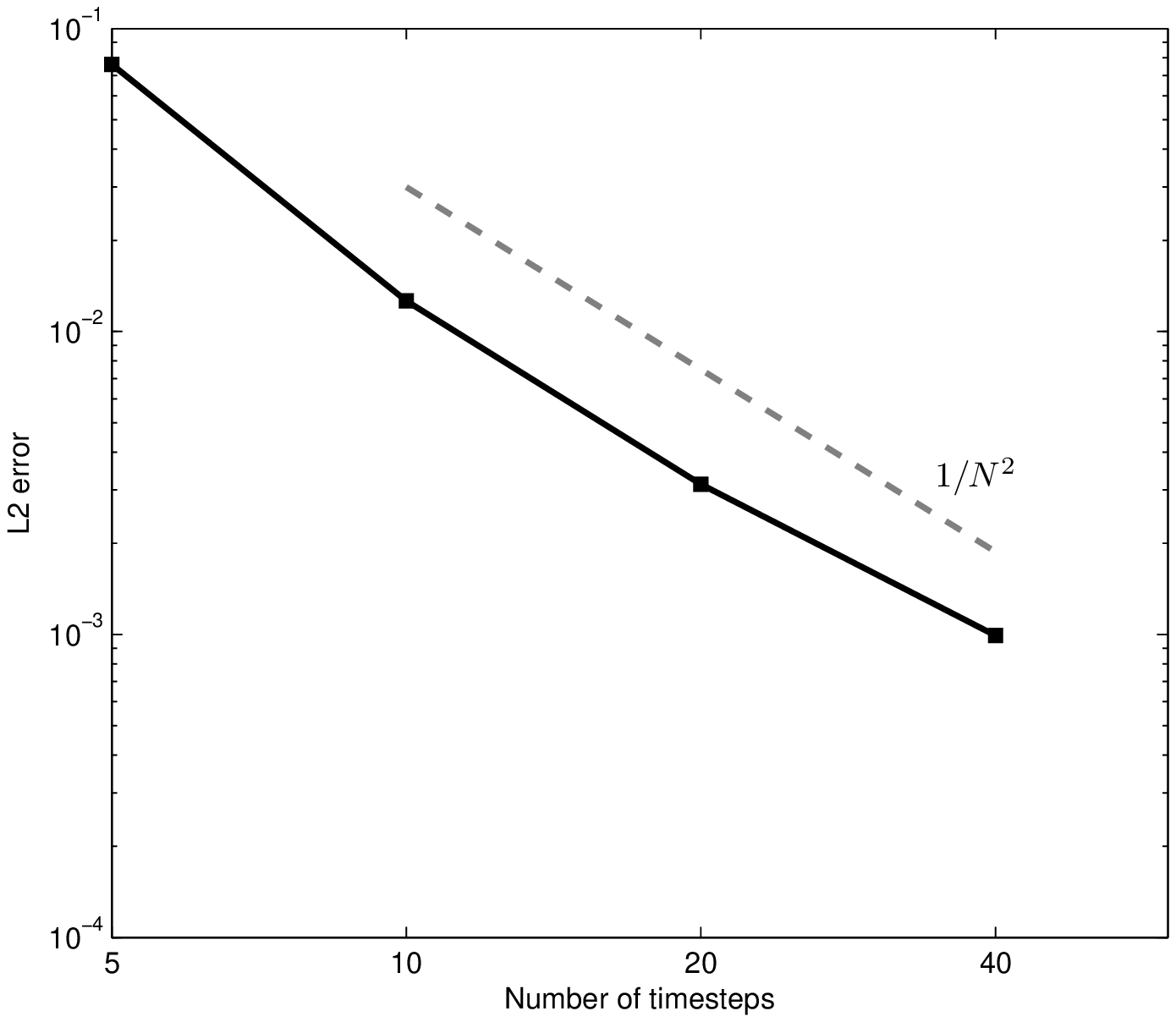}\includegraphics[width=0.5\textwidth]{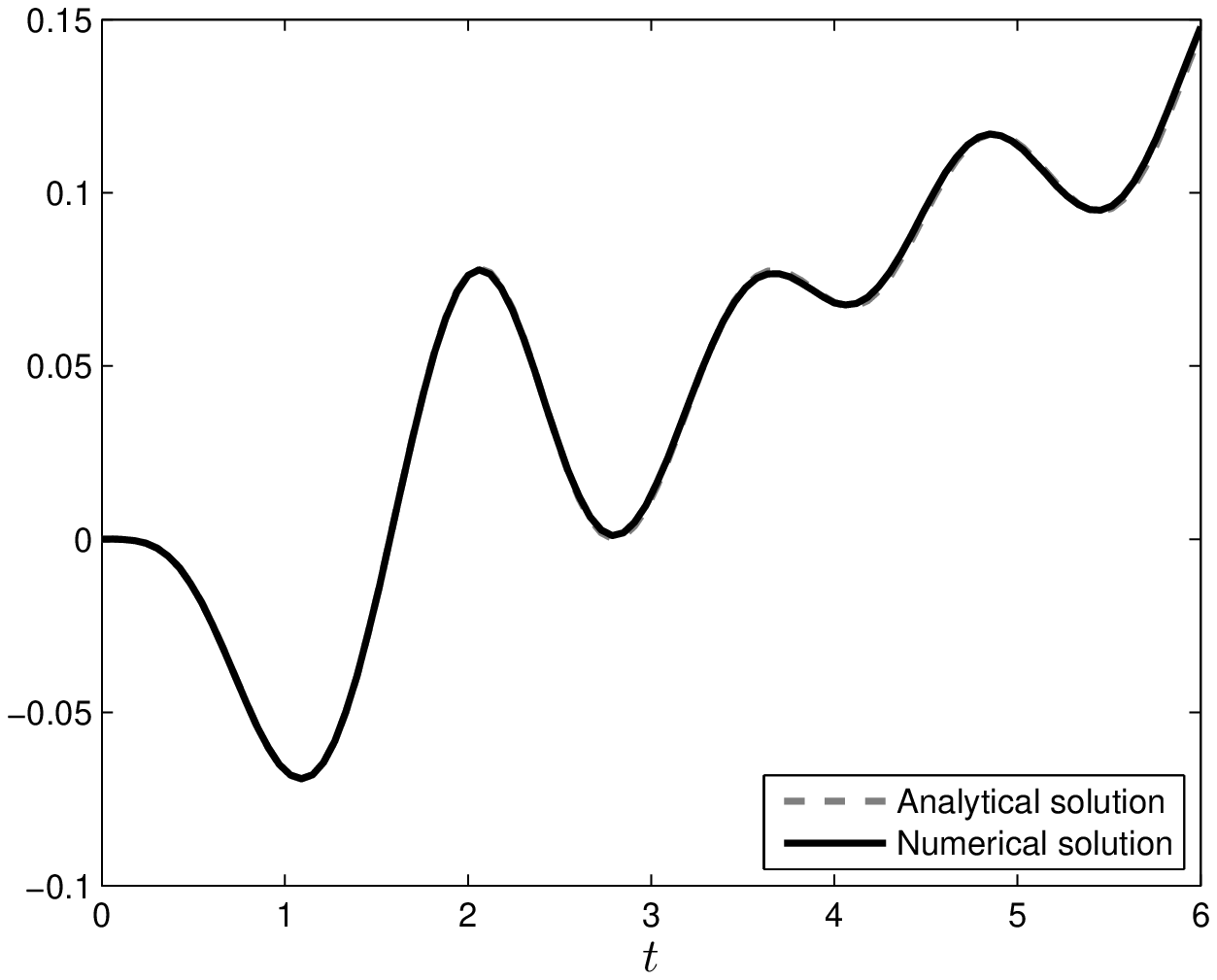}
\end{center}
\caption{Left: log-log scale plot of $\|\Phi-\Phi_{\mathrm{Galerkin}}\|_{L^2(\Gamma;\,[0,T])}$ for $g(t)=\sin(3t)t^2\e^{-t}, n=0, T=6$, and local polynomial approximation space of degree $p=2$. Right: numerical solution in $x\in\Gamma$ obtained using 20 timesteps and $p=2$.}
\label{fig:error2}
\end{figure}

In the following set of numerical experiments we test the convergence of the method using the reference solutions introduced above. First, we solve \eqref{BIE} on the unit sphere $\Sp^2$ for $t\in[0,6]$ with the purely time-dependent right-hand side of the form 
\[g(x,t):=\sin(3t)t^2\e^{-t}Y_0^0(x).\] 
Note that $Y_0^0$ is constant on the unit sphere. The error of the numerical solution is measured in the $\mathrm{L}^2(\Gamma;[0,6])$ norm. Using the local (temporal) polynomial approximation space of order $p=1$ we obtain the convergence rate of $N^{-1}$ (see Figure \ref{fig:error1}); the convergence for $p=2$ is depicted in~Figure~\ref{fig:error2}. 

\begin{figure}[htb]
\begin{center}
\includegraphics[width=0.5\textwidth]{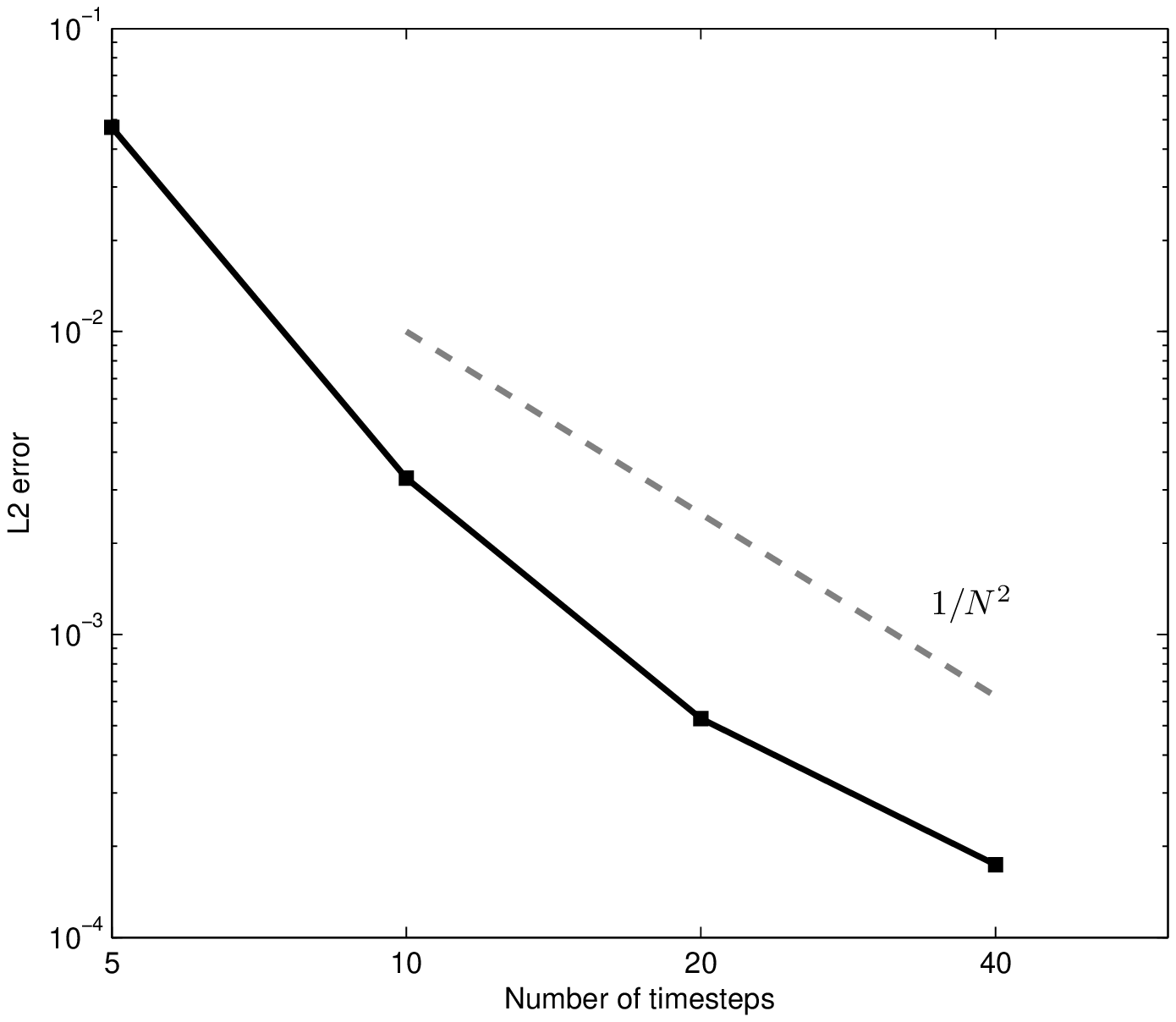}\includegraphics[width=0.5\textwidth]{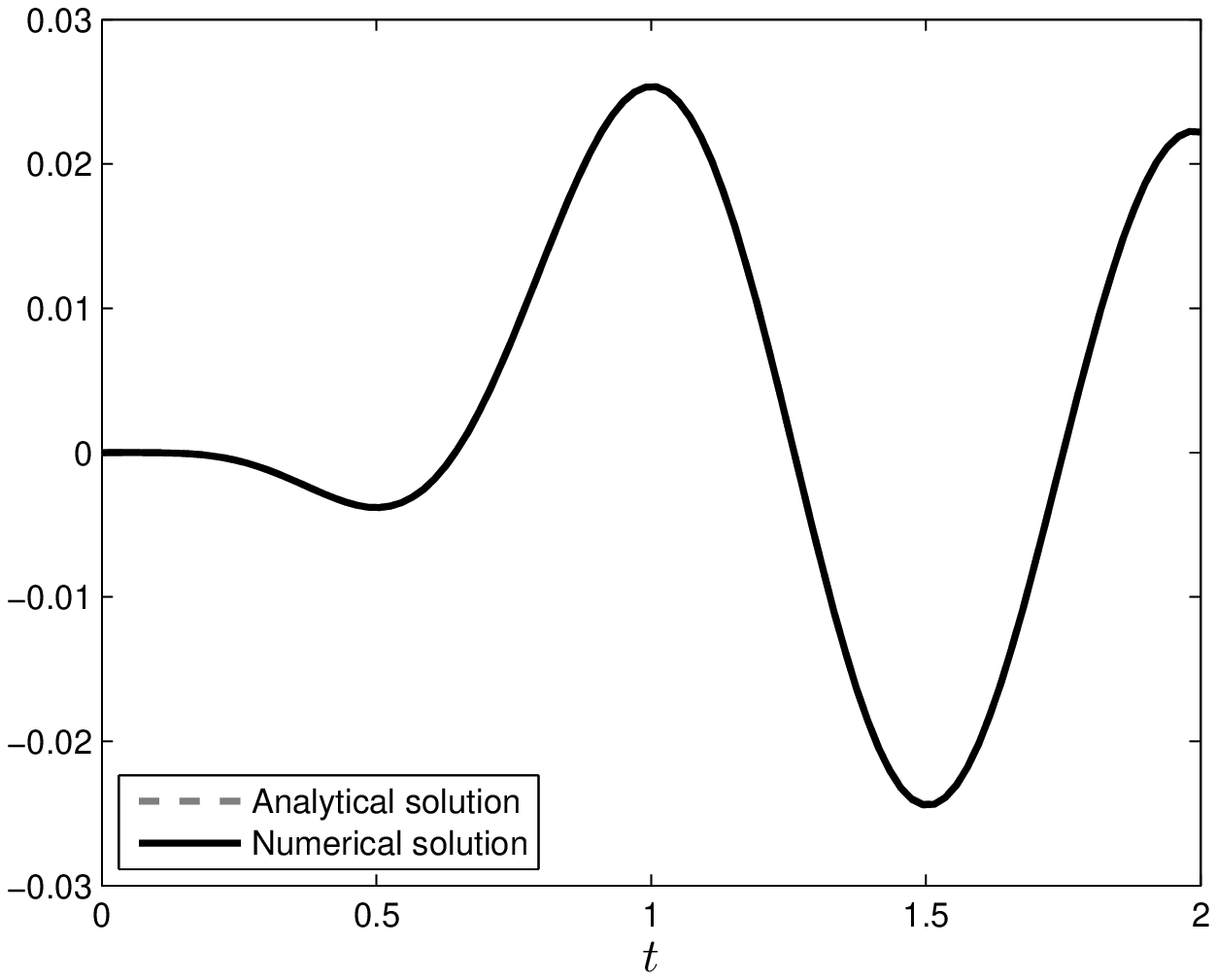}
\end{center}
\caption{Left: log-log scale plot of $\|\Phi-\Phi_{\mathrm{Galerkin}}\|_{L^2(\Gamma;\,[0,T])}$ for $g(t)=\sin(2\pi t)t^3\e^{-2t}, n=1, T=2$, and local polynomial approximation space of degree $p=2$. Right: numerical solution in $x\in\Gamma$ obtained using 20 timesteps and $p=2$.}
\label{fig:error3}
\end{figure}

For the second numerical experiment we consider the right-hand side
\[g(x,t):=\sin(2\pi t)t^3\e^{-2t}Y_1^0(x),\]
and solve \eqref{BIE} with $\Gamma:=\Sp^2$ and $t\in[0,2]$. The numerical solution evaluated in the point $\Gamma\ni x=(0,0,1)$ obtained using the temporal approximation space of order $p=2$ and its convergence to the analytical solution in the $L^2(\Gamma;[0,2])$ norm are depicted in Figure~\ref{fig:error3}.

\subsection{Convergence of iterative solvers}

\begin{figure}[htb]
\begin{center}
\includegraphics[width=0.6\textwidth]{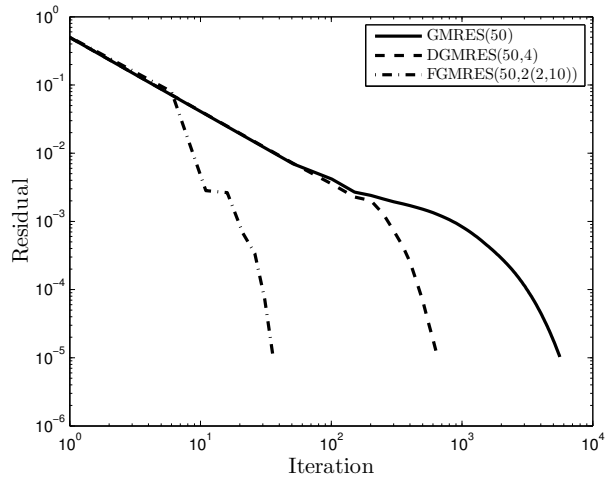}
\end{center}
\caption{Convergence history of GMRES, DGMRES, and FGMRES with recursive preconditioner ($p=1, N=20$).}
\label{fig:dgmres}
\end{figure}

We demonstrate the performance of the iterative solvers introduced in Section~\ref{Sec:SolLinSys} on a set of numerical experiments. The tests were performed using a mesh with 320 surface elements, the relative precision of the solvers is set to $\varepsilon:=10^{-5}$. The results for the local polynomial approximation space of degree $p=1$ are shown in Tables~\ref{tab:gmres1} and \ref{tab:fgmres1}; in Tables~\ref{tab:gmres2} and \ref{tab:fgmres2} we provide results for $p=2$. The convergence history of all solvers is depicted in Figure~\ref{fig:dgmres}. The arguments $m$ and $l$ in GMRES($m$) and DGMRES($m, l$) denote the number of iterations between restarts and the number of vectors added to the basis of the invariant subspace after each restart, respectively. Moreover, by FGMRES($m, r(i_1, \ldots, i_r)$) we denote the solution by flexible GMRES preconditioned as in Section~\ref{Sec:AlgPrecond} with $r$ levels of recursion and $i_1, \ldots, i_r$ iterations of the inner FGMRES solver on each level. 

The experiments demonstrate that both preconditioners presented in Section~\ref{Sec:SolLinSys} lead to a significant reduction of iterations and computational time. Moreover, the results indicate that only a few iterations of the inner solver are necessary during the application of the preconditioner in the case of FGMRES. The best convergence was obtained when the solution of \eqref{innerSystems} was approximated by two iterations of FGMRES on the first level and ten iterations on the second level of recursion. 

%In Table~\ref{tab:gmres1} and \ref{tab:fgmres1} we compare solution times and iteration counts of GMRES, DGMRES, and FGMRES preconditioned as in Section~\ref{Sec:AlgPrecond} with respect to the number of timesteps.  For the test we used the local polynomial approximation space of degree $p=1$. Both DGMRES and FGMRES led to a significant reduction in number of iterations. 

%In Table~\ref{tab:gmres2} we summarize the convergence of GMRES and DGMRES for $p=2$. Again, the number of iterations and solution times were noticeably reduced. Although the usage of prec-GMRES in this case led to a significant reduction in iterations, the overall computational time increased due to costly application of the preconditioner. Therefore, its further modifications are necessary to ensure the applicability to a wider range of problems.

%\begin{table}[htb]
%\begin{center}
%\begin{tabular}{lcccccc}
% \toprule
% &\multicolumn{2}{c}{GMRES($50$)} & \multicolumn{2}{c}{DGMRES($50,4$)} & \multicolumn{2}{c}{prec-GMRES($50$)} \\
%
%$N$ & \# iterations & time [s] & \# iterations & time [s] & \# iterations & time [s] \\
% \midrule
%5 	& \hphantom{0}595 	& \hphantom{0}1.6	& 246 & \hphantom{0}0.8 & 6 & \hphantom{0}2.7 \\
%10 	& 2121 	& 14.3	& 421 & \hphantom{0}3.3 & 6 & \hphantom{0}9.4 \\
%15 	& 4021  & 44.5 	& 580 & \hphantom{0}8.5 & 6 & 24.8 \\
%20 	& 5448  & 99.0	& 510 & 14.9 & 7 & 34.8 \\
%\bottomrule
%\end{tabular}
%\caption{Convergence of iterative solvers for $p=1$.}
%\label{tab:gmres1}
%\end{center}
%\end{table}

\begin{table}[bt]
\begin{center}
\begin{tabular}{lcccc}
 \toprule
 &\multicolumn{2}{c}{GMRES($50$)} & \multicolumn{2}{c}{DGMRES($50,4$)}  \\

$N$ & \# iterations & time [s] & \# iterations & time [s] \\
 \midrule
5 	& \hphantom{0}595 	& \hphantom{0}1.6	& 246 & \hphantom{0}0.8 \\
10 	& 2121 	& 14.3	& 421 & \hphantom{0}3.3 \\
15 	& 4021  & 44.5 	& 580 & \hphantom{0}8.5 \\
20 	& 5448  & 99.0	& 510 & 14.9 \\
\bottomrule
\end{tabular}
\caption{Convergence of GMRES and DGMRES for $p=1$.}
\label{tab:gmres1}
\end{center}
\end{table}

\begin{table}[bht]
\begin{center}
\begin{tabular}{lcccccc}
 \toprule
 &\multicolumn{2}{c}{FGMRES($50,1(10)$)} & \multicolumn{2}{c}{FGMRES($50,1(5)$)} & \multicolumn{2}{c}{FGMRES($50,2(2,10)$)} \\

$N$ & \# iterations & time [s] & \# iterations & time [s] & \# iterations & time [s] \\
 \midrule
5 	& 24 	& 0.7	& \hphantom{0}45 & \hphantom{0}0.9 & 23 & \hphantom{0}0.8 \\
10 	& 43 	& 3.1	& 126 & \hphantom{0}6.8 & 26 & \hphantom{0}3.3 \\
15 	& 51  	& 7.3 	& 205 & 20.0 & 28 & \hphantom{0}5.9\\
20 	& 48  	& 9.7	& 341 & 51.2 & 34 & 10.6 \\
\bottomrule
\end{tabular}
\caption{Convergence of FGMRES with recursive preconditioner for $p=1$.}
\label{tab:fgmres1}
\end{center}
\end{table}

\begin{table}[htb]
\begin{center}
\begin{tabular}{lcccc}
\toprule
 &\multicolumn{2}{c}{GMRES($50$)} & \multicolumn{2}{c}{DGMRES($50,2$)} \\

$N$ & \# iterations & time [s] & \# iterations & time [s] \\
\midrule
5 	& 1985 	& \hphantom{0}14.3	& 1434 & 11.6 \\
10 	& 5404 	& \hphantom{0}89.1	& 3834 & 67.9 \\
15 	& 4634  & 132.7 	& 1491 & 52.1 \\
20 	& 6383  & 293.3	& 1269 & 68.1 \\
\bottomrule
\end{tabular}
\caption{Convergence of GMRES and DGMRES for $p=2$.}
\label{tab:gmres2}
\end{center}
\end{table}

\begin{table}[htb]
\begin{center}
\begin{tabular}{lcccccc}
 \toprule
 &\multicolumn{2}{c}{FGMRES($50,1(10)$)} & \multicolumn{2}{c}{FGMRES($50,1(5)$)} & \multicolumn{2}{c}{FGMRES($50,2(2,10)$)} \\

$N$ & \# iterations & time [s] & \# iterations & time [s] & \# iterations & time [s] \\
 \midrule
5 	& \hphantom{0}83 	& \hphantom{0}5.5	& 210 & \hphantom{00}8.9 & --- & --- \\
10 	& 257 	& 44.3	& 627 & \hphantom{0}77.3 & 94 & 26.6 \\
15 	& 148  	& 48.3 	& 363 & \hphantom{0}82.2 & 56 & 24.8\\
20 	& \hphantom{0}91  	& 45.9	& 399 & 139.4 & 63 & 40.4 \\
\bottomrule
\end{tabular}
\caption{Convergence of FGMRES with recursive preconditioner for $p=2$.}
\label{tab:fgmres2}
\end{center}
\end{table}

\subsection{Wave scattering off a submarine}

\begin{figure}[t]
\begin{center}
\includegraphics[width=0.6\textwidth,natwidth=1351,natheight=387]{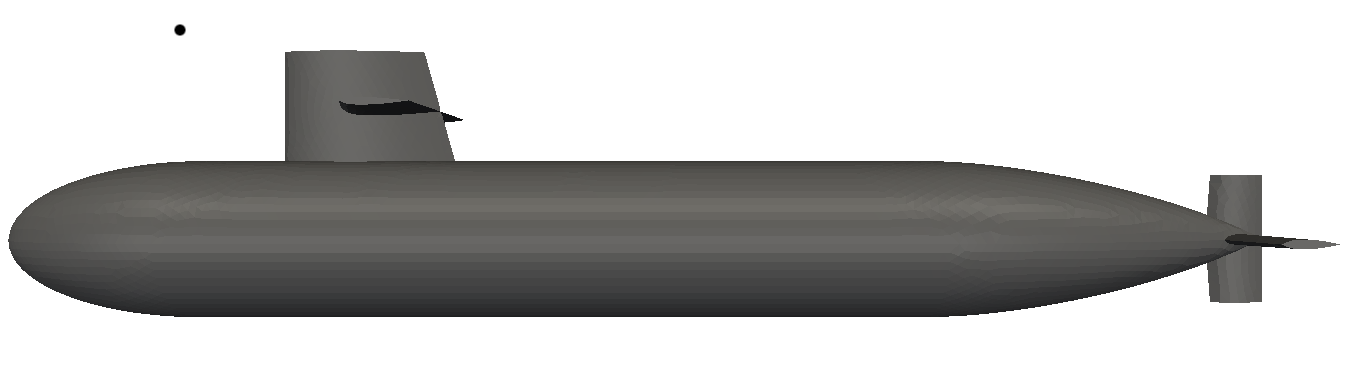}
\end{center}
\caption{Submarine-shaped computational geometry.}
\label{fig:submarine}
\end{figure}

The following numerical experiments serve to demonstrate the behavior of the solver on more complex geometries, such as the submarine-shaped scatterer depicted in Figure~\ref{fig:submarine}. We solve the sound-hard scattering problem for $t\in[0,6]$ with the planar incident wave \cite{Sayas,Glaefke2012}
\[
u^{\mathrm{inc}}\left( x, t\right)  :=\left\{
\begin{array}
[c]{ll}%
A\cos(kx+\varphi_0-\omega t), & \omega t-m_f\geq kx \geq \omega t-m_t,\\
0, & \mathrm{otherwise},\\
\end{array}
\right.
\]
where $A=0.02$, $k=(-\pi/\sqrt{2}, 0.0, -\pi/\sqrt{2})$, $\omega=\pi$, $m_f=6\pi$, and $m_t=8\pi$. The Neumann boundary condition \eqref{BoundaryConditions} is given by
\[\partial_n u=-\frac{\partial u^{\mathrm{inc}}}{\partial n}=
\left\{
\begin{array}
[c]{ll}%
A\sin(kx+\varphi_0-\omega t)kn, & \omega t-m_f\geq kx \geq \omega t-m_t,\\
0, & \mathrm{otherwise},\\
\end{array}
\right.\]
with $n$ being the unit outward normal vector to $\Gamma$. We decompose the boundary of the scatterer into 5604 surface elements, the time interval into 40 equidistant time steps and use the local polynomial approximation space of order $p=1$. 

\begin{figure}[htb]
\begin{center}
\includegraphics[width=0.49\textwidth]{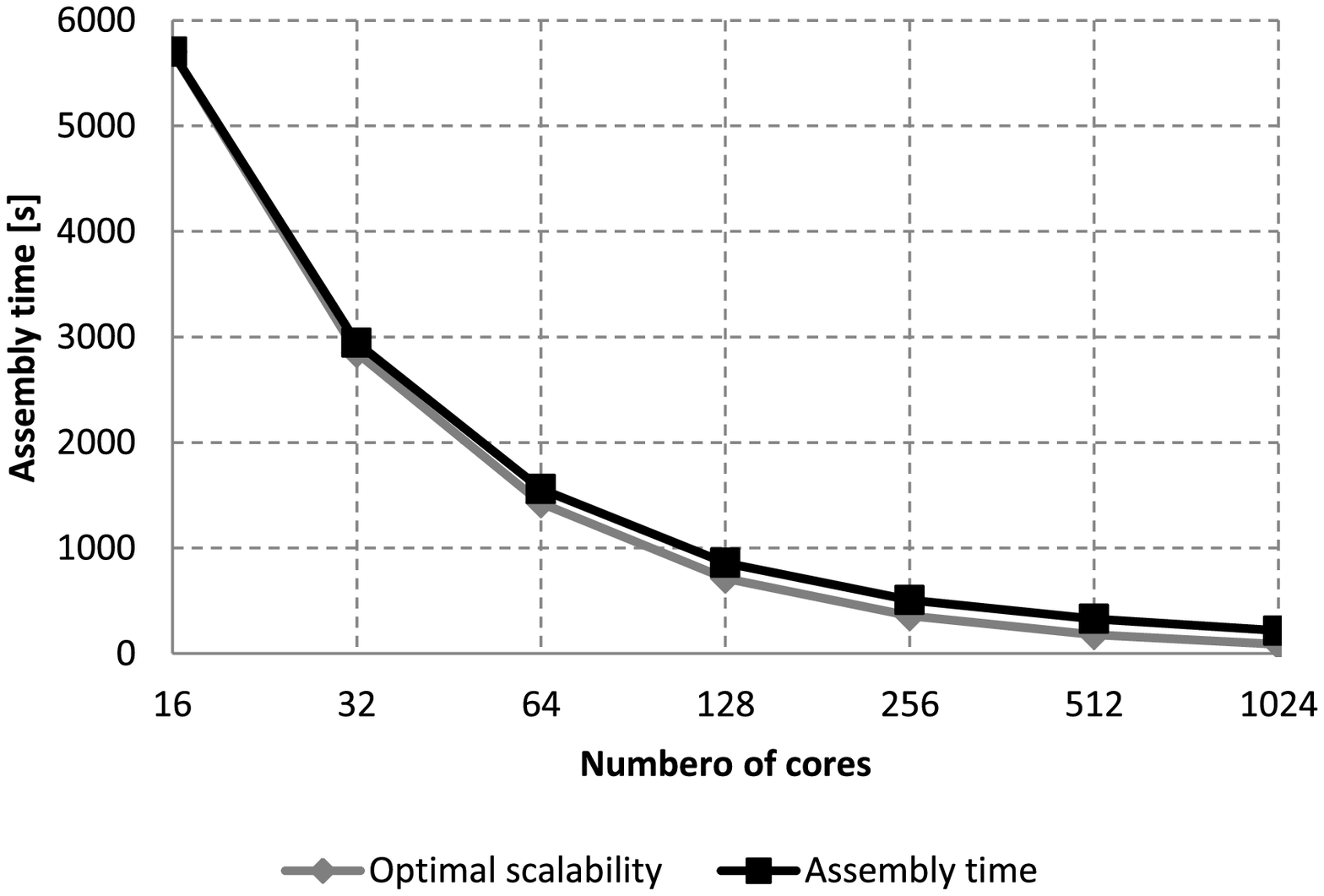}
\includegraphics[width=0.49\textwidth]{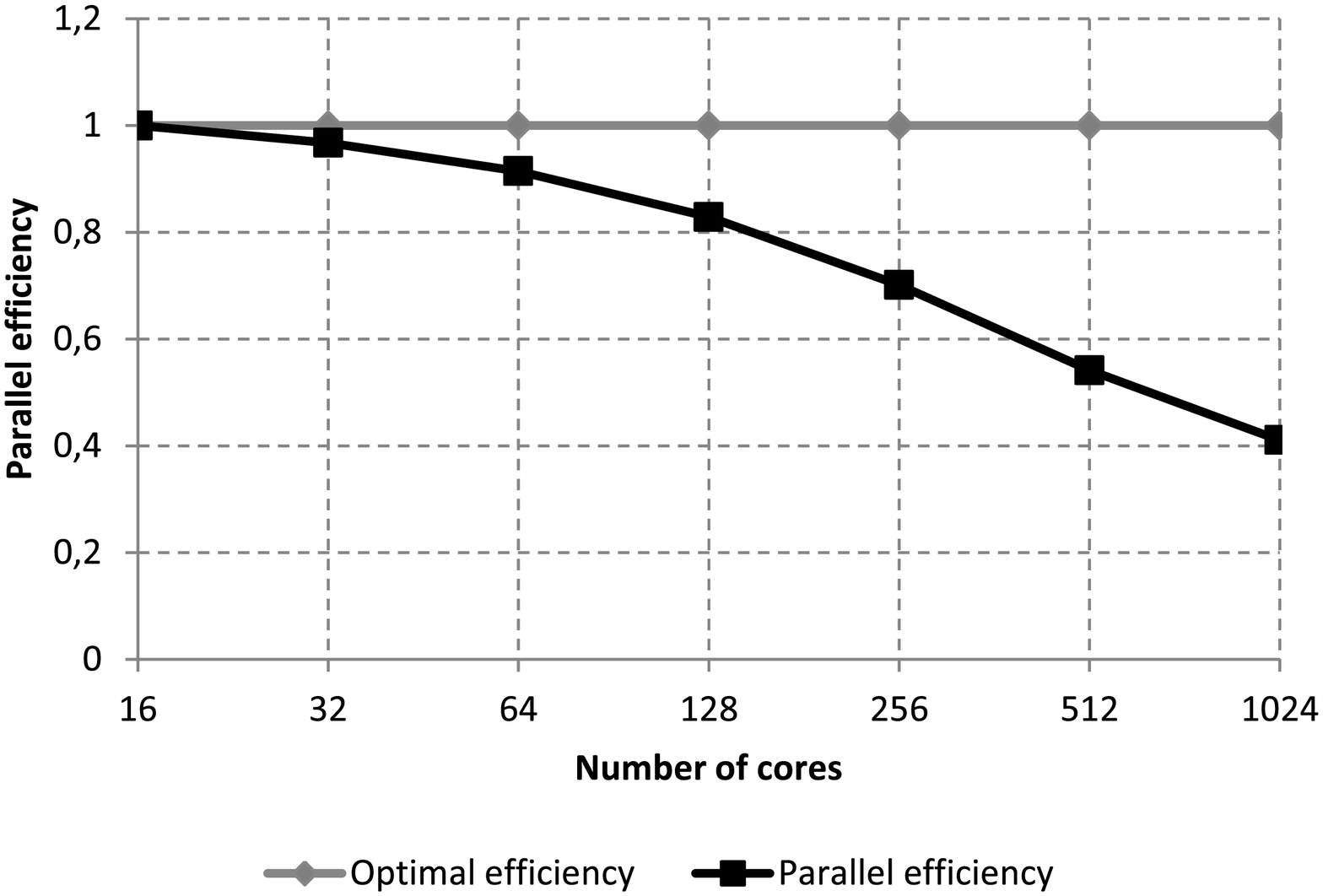}
\end{center}
\caption{Parallel scalability and efficiency of the system matrix assembly up to 1024 cores (64 nodes).}
\label{fig:subm_scal}
\end{figure}

On the left-hand side of Figure \ref{fig:subm_scal} the scalability of the system matrix assembly up to 1024 cores is depicted. The assembly time was reduced from 5702 s on 16 cores to 217 s on 1024 cores. The corresponding parallel efficiency with respect to a single computational node is depicted in the right-hand side of Figure \ref{fig:subm_scal}. The main bottleneck influencing the scalability is currently caused by the MPI communication during the non-zero pattern computation and during the gathering of matrix contributions from MPI processes.

\begin{figure}[htb]
\begin{center}
\includegraphics[width=0.85\textwidth]{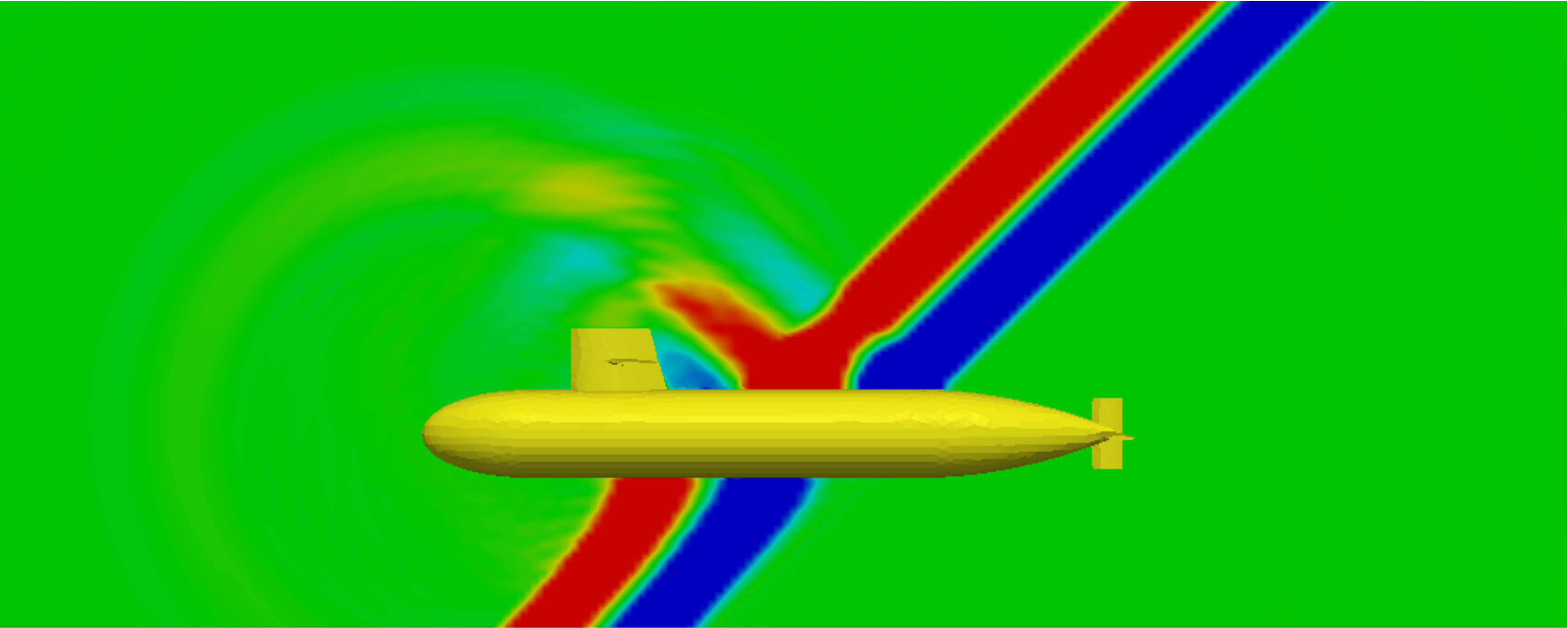}
\end{center}
\caption{Solution (sum of incident and scattered wave) at time $t=3.5\,\mathrm{s}$.}
\label{fig:sol_submarine}
\end{figure}

Finally, Figure \ref{fig:sol_submarine} depicts the solution (sum of incident and scattered wave) in the vicinity of the scatterer at time $t=3.5 \,\mathrm{s}$. The double layer potential was evaluated in 66049 points for 40 timesteps. The evaluation took $327 \,\mathrm{s}$ using 64 computational nodes.

\section{Conclusion}
We considered the numerical modelling of time-dependent, three-dimensional sound-hard scattering phenomena in unbounded domains. We used time-domain boundary integral equations in combination with space-time Galerkin methods to solve the arising problems numerically. The use of $C^\infty$-smooth and compactly supported temporal basis functions simplifies the generation of the system matrix and as a consequence allows the use of higher order approximation schemes in time which significantly improves the accuracy of the method. Due to the overlap of the basis functions in time the arising linear systems that need to be solved admit a block Hessenberg structure. We examined the performance of a restarted GMRES solver preconditioned by deflation  and proposed an algebraic preconditioner that exploits the block Hessenberg structure of the matrix. Various examples show that both solvers require a considerably lower number of iterations than a conventional GMRES method. Especially the second approach is promising since the overall solution time could be significantly reduced. A rigorous analysis and a possible parallel implementation of the preconditioner, however, will be considered in the future.

We furthermore presented a parallel implementation of the scheme that is based on the boundary element library BEM4I. It exploits the special structure of the system matrix to reduce computational complexity and memory requirements. It uses OpenMP and MPI for parallelization in shared and distributed memory. Numerical experiments show good parallel scalability and efficiency of the computations in a distributed memory architecture with up to $1024$ cores.

\section*{Acknowledgement}
The first author gratefully acknowledges the support given by the Swiss National Science Foundation (SNF), No. P2ZHP2\_148705. Latter authors were supported by the IT4Innovations Centre of Excellence project (CZ.1.05/ 1.1.00/02.0070), funded by the European Regional Development Fund and the national budget of the Czech Republic via the Research and Development for Innovations Operational Programme, as well as Czech Ministry of Education, Youth and Sports via the project Large Research, Development and Innovations Infrastructures (LM2011033). The second and third authors were supported by V\v{S}B-TU Ostrava under the grant SGS SP2015/160. The last author was supported by V\v{S}B-TU Ostrava under the grant SGS SP2015/100.

\FloatBarrier

\bibliographystyle{plain}
%\bibliography{mybib}
\bibliography{merta_bib}

\begin{thebibliography}{10}

\bibitem{Babu}
I.~Babu\v{s}ka and J.~Melenk.
\newblock The {P}artition of {U}nity {M}ethod.
\newblock {\em Int. J. Numer. Meths. Eng.}, 40:727--758, 1997.

\bibitem{BamDuongSoft}
A.~Bamberger and T.~{Ha Duong}.
\newblock Formulation variationnelle espace-temps pour le calcul par potentiel
  retard\'{e} de la diffraction d'une onde acoustique ({I}).
\newblock {\em Math. Meth. in the Appl. Sci.}, 8(3):405--435, 1986.

\bibitem{BamDuong}
A.~Bamberger and T.~{Ha Duong}.
\newblock Formulation {V}ariationnelle pour le {C}alcul de la {D}iffraction
  d'une {O}nde {A}coustique par une {S}urface {R}igide.
\newblock {\em Math. Meth. in the Appl. Sci.}, 8(4):598--608, 1986.

\bibitem{BanjaiSchanz}
L.~Banjai and M.~Schanz.
\newblock Wave propagation problems treated with convolution quadrature and
  {BEM}.
\newblock In U.~Langer, M.~Schanz, O.~Steinbach, and W.~L. Wendland, editors,
  {\em Fast Boundary Element Methods in Engineering and Industrial
  Applications}, volume~63 of {\em Lecture Notes in Applied and Computational
  Mechanics}, pages 145--184. Springer Berlin Heidelberg, 2012.

\bibitem{Clenshaw}
I.~Clenshaw.
\newblock A note on the summation of chebyshev series.
\newblock {\em Math. Comp.}, 9:118--120, 1955.

\bibitem{costabel2004time}
M.~Costabel.
\newblock Time-{D}ependent {P}roblems with the {B}oundary {I}ntegral {E}quation
  {M}ethod.
\newblock {\em Encyclopedia of Computational Mechanics}, 2004.

\bibitem{ErhBurPoh96}
J.~Erhel, K.~Burrage, and B.~Pohl.
\newblock Restarted {GMRES} preconditioned by deflation.
\newblock {\em J. Comput. Appl. Math.}, 69(2):303--318, May 1996.

\bibitem{Geranmayeh}
A.~Geranmayeh.
\newblock {\em Time Domain Boundary Integral Equations Analysis}.
\newblock PhD thesis, Fachbereich Elektrotechnik und Informationstechnik,
  Technische Universit¨at Darmstadt, 2011.

\bibitem{Glaefke2012}
M.~Glaefke.
\newblock {\em Adaptive Methods for Time Domain Boundary Integral Equations}.
\newblock PhD thesis, Brunel University, 2012.

\bibitem{HaDuong}
T.~{Ha Duong}.
\newblock On retarded potential boundary integral equations and their
  discretisation.
\newblock In {\em Topics in {C}omputational {W}ave {P}ropagation: {D}irect and
  {I}nverse {P}roblems}, volume 31 of \textit{Lect. Notes Comput. Sci. Eng.},
  pages 301--336. Springer, Berlin, 2003.

\bibitem{DuongLudwigTerrasse}
T.~Ha-Duong, B.~Ludwig, and I.~Terrasse.
\newblock A galerkin bem for transient acoustic scattering by an absorbing
  obstacle.
\newblock {\em International Journal for Numerical Methods in Engineering},
  57(13):1845--1882, 2003.

\bibitem{huang1989some}
Y.~Huang and H.~van~der Vorst.
\newblock {\em Some {O}bservations on the {C}onvergence {B}ehavior of {GMRES}}.
\newblock Reports of the Faculty of Mathematics and Informatics. Delft
  University of Technology. Delft University of Technology, 1989.

\bibitem{KhSaVe}
B.~N. Khoromskij, S.~Sauter, and A.~Veit.
\newblock Fast {Q}uadrature {T}echniques for {R}etarded {P}otentials {B}ased on
  {TT}/{QTT} {T}ensor {A}pproximation.
\newblock {\em Comp. Meth. Appl. Math.}, 11(3), 2011.

\bibitem{BEM4I}
M.~Merta and J.~Zapletal.
\newblock {BEM4I Library}.
\newblock \url{http://industry.it4i.cz/en/products/bem4i/}.
\newblock Accessed: 2015-02-26.

\bibitem{nedelec2001acoustic}
J.~N\'{e}d\'{e}lec.
\newblock {\em Acoustic and {E}lectromagnetic {E}quations: {I}ntegral
  {R}epresentations for {H}armonic {P}roblems}.
\newblock Number 144 in Acoustic and electromagnetic equations: integral
  representations for harmonic problems. Springer, 2001.

\bibitem{saad1993}
Y.~Saad.
\newblock A {F}lexible {I}nner-outer {P}reconditioned {GMRES} {A}lgorithm.
\newblock {\em SIAM J. Sci. Comput.}, 14(2):461--469, Mar. 1993.

\bibitem{Saad1986gmres}
Y.~Saad and M.~Schultz.
\newblock {GMRES}: A generalized minimal residual algorithm for solving
  nonsymmetric linear systems.
\newblock {\em SIAM Journal on scientific and statistical computing},
  7(3):856--869, 1986.

\bibitem{Sauter5}
S.~Sauter and A.~Veit.
\newblock A {G}alerkin {M}ethod for {R}etarded {B}oundary {I}ntegral
  {E}quations with {S}mooth and {C}ompactly {S}upported {T}emporal {B}asis
  {F}unctions.
\newblock {\em Numerische Mathematik}, pages 1--32, 2012.

\bibitem{SauterVeit2013}
S.~Sauter and A.~Veit.
\newblock Adaptive time discretization for retarded potentials.
\newblock {\em Preprint 06-2013, University of Zurich}, 2013.

\bibitem{SauterVeitIMA}
S.~Sauter and A.~Veit.
\newblock Retarded {B}oundary {I}ntegral {E}quations on the {S}phere: {E}xact
  and {N}umerical {S}olution.
\newblock {\em IMA J. Numer. Anal.}, 2(34):675--699, 2014.

\bibitem{Sayas}
F.~Sayas.
\newblock Retarded {P}otentials and {T}ime {D}omain {B}oundary {I}ntegral
  {E}quations: a road-map.
\newblock Lecture notes, 2013.

\bibitem{Veit_diss}
A.~Veit.
\newblock {\em Numerical {M}ethods for {T}ime-{D}omain {B}oundary {I}ntegral
  {E}quations}.
\newblock PhD thesis, University of Zurich, 2012.

\end{thebibliography}

\end{document}